\input amstex\input epsf\input lpic-mod
\documentstyle{amsppt}\nologo\footline={}\subjclassyear{2000}

\font\oitossi=cmssqi8
\font\oitoss=cmssq8

\def\B{\mathop{\text{\rm B}}}
\def\S{\mathop{\text{\rm S}}}
\def\E{\mathop{\text{\rm E}}}
\def\PU{\mathop{\text{\rm PU}}}
\def\L{\mathop{\text{\rm L}}}
\def\G{\mathop{\text{\rm G}}}
\def\dist{\mathop{\text{\rm dist}}}
\def\SU{\mathop{\text{\rm SU}}}
\def\Stab{\mathop{\text{\rm Stab}}}
\def\U{\mathop{\text{\rm U}}}
\def\su{\mathop{\text{\rm su}}}
\def\T{\mathop{\text{\rm T}}}
\def\Lin{\mathop{\text{\rm Lin}}}
\def\tr{\mathop{\text{\rm tr}}}
\def\ta{\mathop{\text{\rm ta}}}
\def\Im{\mathop{\text{\rm Im}}}
\def\Re{\mathop{\text{\rm Re}}}
\def\PGL{\mathop{\text{\rm PGL}}}

\hsize450pt\topmatter\title Reflections, bendings, and
pentagons\endtitle\author Sasha Anan$'$in\endauthor\address
Departamento de Matem\'atica, IMECC, Universidade Estadual de
Campinas,\newline13083-970--Campinas--SP, Brasil\endaddress\subjclass
30F60 (51M10, 57S30)\endsubjclass\abstract We study relations between
reflections in (positive or negative) points in the complex hyperbolic
plane. It is easy to see that the reflections in the points $q_1,q_2$
obtained from $p_1,p_2$ by moving $p_1,p_2$ along the geodesic
generated by $p_1,p_2$ and keeping the (dis)tance between $p_1,p_2$
satisfy the {\it bending relation\/} $R(q_2)R(q_1)=R(p_2)R(p_1)$. We
show that a generic isometry $F\in\SU(2,1)$ is a product of $3$
reflections, $F=R(p_3)R(p_2)R(p_1)$, and describe all such
decompositions: two decompositions are connected by finitely many
bendings involving $p_1,p_2$/$p_2,p_3$ and geometrically equal
decompositions differ by an isometry centralizing $F$.

Any relation between reflections gives rise to a representation
$H_n\to\PU(2,1)$ of the {\it hyperelliptic group\/} $H_n$ generated by
$r_1,\dots,r_n$ with the defining relations $r_n\dots r_1=1$,
$r_j^2=1$. The theorem mentioned above is essential to the study of the
Teichm\"uller space $\Cal TH_n$. We describe all nontrivial
representations of $H_5$, called {\it pentagons,} and conjecture
that they are faithful and discrete.\endabstract\endtopmatter\document

\rightline{\oitossi What do you feed the people generated by
reflections?}

\medskip

\rightline{\oitoss --- A.~A.~Kirillov to E.~B.~Vinberg on the way to
Baikal Lake}

\bigskip

\bigskip

\centerline{\bf1.~Introduction}

\medskip

We will deal with the complex hyperbolic plane $\Bbb H_\Bbb C^2$. In
our notation, we follow [AGr] and~[AGG]. Let $V$ denote a
$3$-dimensional $\Bbb C$-linear space equipped with a hermitian form
$\langle-,-\rangle$ of signature $++-$. By definition,
$\Bbb H_\Bbb C^2:=\B V$ and $\partial\Bbb H_\Bbb C^2:=\S V$, where
$$\B V:=\big\{p\in\Bbb P_\Bbb CV\mid\langle p,p\rangle<0\big\},\qquad\S
V:=\big\{p\in\Bbb P_\Bbb CV\mid\langle p,p\rangle=0\big\},\qquad\E
V:=\big\{p\in\Bbb P_\Bbb CV\mid\langle p,p\rangle>0\big\}$$
(see [AGr] and [AGG] for details).

\smallskip

Let us start with some heuristic considerations. For many geometries,
it is well known that every isometry is a product of involutions. It is
curious to understand to what extent such a decomposition is unique. In
other words, we study the `defining relations' $R_n\dots R_2R_1=1$ and
$R_j^2=1$.

Consider $\PU V$, the group of holomorphic isometries of
$\Bbb H_\Bbb C^2$. Every involution in $\PU V$ is a reflection $R(p)$
in some nonisotropic point $p\in\Bbb P_\Bbb CV\setminus\S V$ (positive
or negative).

\smallskip

$\bullet$ The relation $R(p)R(p)=1$ is said to be a {\it cancellation.}

\smallskip

\noindent
Clearly, these are the only relations with $n=2$.

\smallskip

$\bullet$ For pairwise orthogonal
$p_1,p_2,p_3\in\Bbb P_\Bbb CV\setminus\S V$, the reflections $R(p_1)$,
$R(p_2)$, $R(p_3)$ commute and $R(p_3)R(p_2)R(p_1)=1$. We call these
relations {\it orthogonal.}

\smallskip

\noindent
It is easy to see that these are the only relations with $n=3$ (see
Remark 3.2).

\smallskip

Suppose that distinct points $p_1,p_2\in\Bbb P_\Bbb CV\setminus\S V$
generate a hyperbolic\footnote{$\L(p_1,p_2)$ is {\it hyperbolic,
spherical, euclidean\/} if the signature of $\Bbb Cp_1+\Bbb Cp_2$ is
respectively $+-$, $++$, $+0$.}
projective line $\L(p_1,p_2)$. Then the polar point $p$ to
$\L(p_1,p_2)$ is positive and
$R(p_2)R(p_1)=R(p_2)R(p)R(p)R(p_1)=R(p'_2)R(p'_1)$, where
$p'_1,p'_2\in\L(p_1,p_2)$ are respectively orthogonal to $p_1,p_2$. So,
with the use of cancellation and orthogonal relations, we can
simultaneously alter the signs of $p_1,p_2$ if the projective line
$\L(p_1,p_2)$ is hyperbolic.

Given distinct and nonorthogonal points
$p_1,p_2\in\Bbb P_\Bbb CV\setminus\S V$, we can move $p_1,p_2$ along
the geodesic $\G{\wr}p_1,p_2{\wr}$ and keep constant the distance
(tance) between them. The new points $q_1,q_2$ satisfy the relation
$R(q_2)R(q_1)=R(p_2)R(p_1)$. More formally: There exists a
one-parameter subgroup $B:(\Bbb R,+)\to\PU V$ such that $B(s)$
centralizes $R(p_2)R(p_1)$,
$\dist\big(p_1,B(s)p_1\big)=\dist(p_1,p_2)|s|$, $B(1)p_1=p_2$ if
$p_1,p_2$ have the same sign, and $B(1)p_1\in\L(p_1,p_2)$ is orthogonal
to $p_2$ if $p_1,p_2$ have different signs. For $q_1:=B(s)p_1$ and
$q_2:=B(s)p_2$, we have
$R(q_2)R(q_1)=R\big(B(s)p_2\big)R\big(B(s)p_1\big)=
B(s)R(p_2)R(p_1)B(-s)=R(p_2)R(p_1)$.
We~say that the points $q_1,q_2$ are obtained from $p_1,p_2$ by means
of {\it bending.}

\smallskip

$\bullet$ The relation $R(p_2)R(p_1)=R(q_2)R(q_1)$ is called a {\it
bending relation.}

\smallskip

\noindent
It is easy to show that, modulo cancellations and orthogonal relations,
the bending relations provide all relations with $n=4$ (Corollary 3.3).

\smallskip

Dealing with a general relation $R(p_n)\dots R(p_2)R(p_1)=1$, we can
bend $p_{j-1},p_j$, then bend $p_{k-1},p_k$, etc.~(the indices are
modulo $n$). At some moment, we may arrive at $p_{j-1}=p_j$ or
$\langle p_{j-1},p_j\rangle=0$ and cancel $R(p_j)R(p_{j-1})$ or apply
an orthogonal relation, thus diminishing $n$. Also, with the use of
cancellation and orthogonal relations, we can simultaneously alter the
signs of $p_{j-1},p_j$ if the projective line $\L(p_{j-1},p_j)$ is
hyperbolic. We expect that there is a finite number of {\it basic\/}
relations that generate all the relations in the described way.

\smallskip

Denote by $H_n$ the group generated by $r_1,r_2,\dots,r_n$ with the
defining relations $r_j^2=1$, $j=1,\dots,n$, and $r_n\dots r_2r_1=1$.
Every relation $R(p_n)\dots R(p_2)R(p_1)=1$ gives rise to a
representation $H_n\to\PU V$. The most interesting are the faithful and
discrete ones.

\medskip

{\bf1.1.~Conjecture.} There exists a finite number of basic relations,
each one induces a faithful and discrete representation that remains
faithful and discrete after a finite number of bendings.

\medskip

For such a representation, the projective line $\L(p_{j-1},p_j)$ cannot
be spherical as $R(p_j)R(p_{j-1})$ would otherwise be elliptic. Nor
does it rate to be euclidean since bending $p_j,p_{j+1}$ usually leads
to a spherical $\L(p_{j-1},p_j)$. Unfortunately, even if at the
beginning all projective lines $\L(p_{j-1},p_j)$ are hyperbolic, some
of them can become spherical after a few bendings (see Section 4.3).
This will not happen if at most one of the $p_j$'s is positive.

\smallskip

As a first step, we study the representations
$$\varrho:H_5\to\PU V,\qquad\varrho:r_j\mapsto R(p_j),\qquad p_j\in\Bbb
P_\Bbb CV\setminus\S V$$
such that at most one of the $p_j$'s is positive. For brevity, we call
them {\it pentagons.} Lifting a pentagon to $\SU V$, we obtain the
relation
$$R(p_5)R(p_4)R(p_3)R(p_2)R(p_1)=\delta,$$
where $\delta^3=1$. We show that, modulo conjugation, two pentagons
with the same $\delta$ are connected by means of a finite number of
bendings (see Subsection 5.4). In other words, each cubic root of unity
$\delta$ provides a connected component of the space of pentagons. If
$\delta=1$, the pentagons are $\Bbb R$-fuchsian and the $p_j$'s are
negative. In this case, all pentagons are faithful and discrete [ABG].
Therefore, while deforming such a pentagon, there is always some real
plane stable under the action of $\varrho H_5$ and $\varrho$ remains
faithful and discrete. An example of a faithful and discrete pentagon
with $\delta\ne1$ is constructed in [Ana].

\medskip

{\bf1.2.~Conjecture.} Every pentagon is faithful and discrete.

\medskip

We believe that the situation with pentagons is quite similar to that
in Toledo's rigidity theorem (see~[Tol] and [Gol1]). This analogy would
be deeper if one could find a holomorphic section of the example in
[Ana]. In this case, the corresponding holomorphic disc in
$\Bbb H_\Bbb C^2$ would be stable under the action of $H_5$.

\medskip

We can describe bendings in the following more general terms. Let
$p\in\Bbb P_\Bbb CV\setminus\S V$. The bundle $\PU V\to\PU V/\Stab p$,
where $\PU V/\Stab p=\B V$ if $p$ is negative and $\PU V/\Stab p=\E V$
if $p$ is positive, is endowed with the (Cartan) connection whose
horizontal subspaces are orthogonal to cosets with respect to the
Killing form. In other words, any piecewise smooth path
$c:[a,b]\to\Bbb P_\Bbb CV\setminus\S V$ uniquely lifts to a path
$F_c:[a,b]\to\PU V$ such that $F_c(s)c(a)=c(s)$ and $F_c(a)=1$. We call
$F_c$ the {\it isometry following\/} $c$. The concept of path-following
isometry can be illustrated with a snail that walks along the path $c$.
During the walk, the ambient geometry may change from the snail
viewpoint and the snail shell may not suit anymore a new geometry.
Nevertheless, the snail takes no risk to break its shell even if the
shell is as large as the whole space $\Bbb H_\Bbb C^2=\B V$ or $\E V$
(see Proposition 2.4 for the list of basic properties of $F_c$). It is
worthwhile mentioning that this concept is applicable to any classic
geometry (see [AGr]). The isometry following a (suitably parameterized)
geodesic $\G{\wr}p_1,p_2{\wr}$ is exactly the bending $B(s)$ involving
the points $p_1,p_2$ (see Proposition 2.6).

\smallskip

The main result of this paper is Theorem 4.2.2. Briefly speaking, it
claims that a generic (i.e., regular in the sense of Definition 3.6)
isometry $F\in\SU V$ is a product of $3$ reflections,
$F=R(p_3)R(p_2)R(p_1)$, that two such decompositions are connected by
means of finitely many bendings involving $p_1,p_2$/$p_2,p_3$, and that
geometrically equal decompositions differ by an isometry centralizing
$F$. (The case when $p_1,p_2,p_3$ are in a same real plane is simpler
and slightly different.) A certain elementary exposition [ABG] of a
part of the classic Teichm\"uller space theory is intended to be
applicable to the complex hyperbolic plane. In this sense, Theorem
4.2.2 is a first move towards studying the Teichm\"uller space
$\Cal TH_n$ of the hyperelliptic group $H_n$ for the complex hyperbolic
plane. Following the line of the proof of Theorem 4.2.2, it seems
possible to generalize it for $>3$ reflections. Moreover, there may
exist a variant of Theorem 4.2.2 valid for an arbitrary classic
geometry [AGr].

\smallskip

In subsequent articles, we plan to study the Teichm\"uller space
$\Cal TH_n$ in detail. In particular, the~path-following isometry is
intended to provide Toledo-like discrete invariants of discrete
representations and to distinguish some of the examples constructed in
[AGG] that have the same topology and the same Toledo invariant.

\bigskip

\centerline{\bf2.~Path-following isometries}

\medskip

In this section, we introduce path-following isometries and describe
bendings, i.e., isometries following geodesics.

\smallskip

We use the conventions and notation from [AGr] and [AGG]. Depending on
the context, the elements in $V$ may denote points in $\Bbb P_\Bbb CV$.
For the sign of a point $p\in\Bbb P_\Bbb CV$, we use the notation
$\sigma p\in\{-1,0,1\}$. In what follows, $\dot c(s)$ stands for the
derivative of $c(s)$ with respect to $s$. We hope that the same
notation for tangent vectors to curves in $V$ and in $\Bbb P_\Bbb CV$
will produce no confusion since these tangent vectors `live' in totally
different places.

\medskip

{\bf2.1.~Remark.} {\sl Let\/ $c:[a,b]\to\Bbb P_\Bbb CV\setminus\S V$ be
a piecewise smooth path. Then there exists a piecewise smooth lift\/
$c_0:[a,b]\to V$ such that\/
$\big\langle c_0(s),c_0(s)\big\rangle=\sigma\in\{-1,+1\}$ and\/
$\big\langle c_0(s),\dot c_0(s)\big\rangle=0$ for all\/ $s\in[a,b]$.
Such a lift is unique modulo multiplication by a unitary constant\/
$u\in\Bbb C$, $|u|=1$.}

\medskip

{\bf Proof.} It is easy to find a piecewise smooth lift with
$\big\langle c_0(s),c_0(s)\big\rangle=\sigma\in\{-1,+1\}$. It follows
that $\big\langle c_0(s),\dot c_0(s)\big\rangle=if(s)$ for all
$s\in[a,b]$, where $f:[a,b]\to\Bbb R$ is piecewise continuous. We can
find a piecewise smooth function $\varphi:[a,b]\to\Bbb R$ such that
$\dot\varphi(s)\sigma=f(s)$ for all $s\in[a,b]$. Now
$$\Big\langle e^{i\varphi(s)}c_0(s),\big(e^{i\varphi(s)}c_0(s)\big)'
\Big\rangle=\big\langle c_0(s),\dot c_0(s)\big\rangle+\big\langle
e^{i\varphi(s)}c_0(s),e^{i\varphi(s)}i\dot\varphi(s)c_0(s)\big\rangle
=$$
$$=i\Big(f(s)-\dot\varphi(s)\big\langle
c_0(s),c_0(s)\big\rangle\Big)=0.$$

For another lift $g(s)c_0(s)$, $g:[a,b]\to\Bbb C$, it follows from
$\big\langle g(s)c_0(s),g(s)c_0(s)\big\rangle=\sigma$ that
$\big|g(s)\big|=1$ and, from
$\Big\langle g(s)c_0(s),\big(g(s)c_0(s)\big)'\Big\rangle=0$, we deduce
$g(s)\dot g(s)\sigma=0$ and, hence, $\dot g(s)=0$
$_\blacksquare$

\medskip

The lift $c_0$ in Remark 2.1 is said to be {\it normalized.}

\smallskip

Note that a normalized lift of a closed path can be nonclosed. For
instance, for orthonormal positive $p,p'\in V$, put
$c_0(s):=p\cos s+p'\sin s$, $s\in[0,\pi]$ (see Proposition 2.6).
Clearly, $\big\langle c_0(s),c_0(s)\big\rangle=1$,
$\dot c_0(s)=-p\sin s+p'\cos s$, and
$\big\langle c_0(s),\dot c_0(s)\big\rangle=0$. So, $c$ is a closed
path, whereas $c_0(\pi)=-c_0(0)$.

\smallskip

In fact, we obtain a $\U(1)$-connection in the tautological line bundle
over $\Bbb P_\Bbb CV\setminus\S V$.

\medskip

{\bf2.2.~Remark.} {\sl Let\/ $c:[a,b]\to\Bbb P_\Bbb CV\setminus\S V$ be
a piecewise smooth path and let\/ $c_0:[a,b]\to V$ be a normalized lift
of\/ $c$. Then\/
$\dot c(t)=\big\langle-,c_0(t)\big\rangle\sigma\dot c_0(t)$ by\/
{\rm[AGG, Lemma 4.1.4]}
$_\blacksquare$}

\medskip

The group $\SU V$ is given by the `equations' $\det X=1$ and
$\langle Xv_1,Xv_2\rangle=\langle v_1,v_2\rangle$, $v_1,v_2\in V$. As
the Lie algebra $\su V$ is the tangent space $\T_1\SU V$, it is given
inside $\Lin_{\Bbb C}(V,V)$ by the `equations' $\tr Y=0$ and
$\langle Yv_1,v_2\rangle+\langle v_1,Yv_2\rangle=0$, $v_1,v_2\in V$,
which can be also written as $\tr Y=0$ and $Y+Y^*=0$. Identifying
$\Lin_{\Bbb C}(V,V)\simeq V^*\otimes_{\Bbb C}V$, we have
$Y^*=\langle-,v\rangle p$ for $Y:=\langle-,p\rangle v$. So,
$\tr Y=\tr Y^*$ if $\langle v,p\rangle\in\Bbb R$. Clearly,
$\T_F\SU V=F\cdot\su V=\su V\cdot F$ for every $F\in\SU V$.

\medskip

{\bf2.3.~Definition.} For $p\in\Bbb P_\Bbb CV\setminus\S V$ and
$t\in\T_p\Bbb P_\Bbb CV\subset\Lin_{\Bbb C}(V,V)$, denote
$\widehat t:=t-t^*\in\su V$. Given a piecewise smooth path
$c:[a,b]\to\Bbb P_\Bbb CV\setminus\S V$ and some $F_0\in\SU V$, we have
a piecewise continuous path $\widehat{\dot c}:[a,b]\to\su V$. The
Cauchy problem $\dot F=\widehat{\dot c}\cdot F$, $F(a)=F_0$, has a
unique piecewise smooth solution $F:[a,b]\to\SU V$. We call $F$ the
{\it isometry following the path\/} $c$ and {\it beginning with\/}
$F_0$. Denote by $F_c$ the isometry that follows $c$ and begins with
$1$.

\medskip

Denote by $c\cup c':[a,d]\to X$ the concatenation of paths
$c:[a,b]\to X$ and $c':[b,d]\to X$ such that $c(b)=c'(b)$.

\medskip

{\bf2.4.~Proposition.} {\sl Let\/
$c:[a,b]\to\Bbb P_\Bbb CV\setminus\S V$ and\/
$c':[b,d]\to\Bbb P_\Bbb CV\setminus\S V$ be piecewise smooth paths such
that\/ $c(b)=c'(b)$, let\/ $F_0\in\SU V$, and let\/
$\varphi,\psi:[a',b']\to[a,b]$ be piecewise smooth reparameterizations
of\/ $c$ such that\/ $\varphi(a')=a$, $\varphi(b')=b$ and\/
$\psi(a')=b$, $\psi(b')=a$. Denote by\/ $F$ the isometry that follows\/
$c$ and begins with\/ $F_0$. Then

\smallskip

$\bullet$ $F(s)=F_c(s)F_0$ for all\/ $s\in[a,b]$,

$\bullet$ $F_{c\cup c'}=F_c\cup F_{c'}F_c(b)$,

$\bullet$ $F_c\circ\varphi=F_{c\circ\varphi}$,

$\bullet$ $F_c\circ\psi=F_{c\circ\psi}F_c(b)$,

$\bullet$ $F_c(s)\big(c_0(a)\big)=c_0(s)$ for all\/ $s\in[a,b]$,
where\/ $c_0:[a,b]\to V$ stands for a normalized lift of $c$.}

\medskip

{\bf Proof.} We will prove only the last two statements. For
$\alpha:=c\circ\psi$, we have
$\dot\alpha(s)=\dot c\big(\psi(s)\big)\dot\psi(s)$ and
$\widehat{\dot\alpha}(s)=\widehat{\dot c}\big(\psi(s)\big)\dot\psi(s)$.
Since
$(F_c\circ\psi)'(s)=\dot F_c\big(\psi(s)\big)\dot\psi(s)=\widehat{\dot
c}\big(\psi(s)\big)F_c\big(\psi(s)\big)\dot\psi(s)=
\widehat{\dot\alpha}(s)(F_c\circ\psi)(s)$
and $(F_c\circ\psi)(a')=F_c(b)$, the fact follows.

\smallskip

For $\beta(s):=F_c(s)\big(c_0(a)\big)$, we have
$$\dot\beta(s)=\dot F_c(s)\big(c_0(a)\big)=\widehat{\dot
c}(s)F_c(s)\big(c_0(a)\big)=\Big(\dot c(s)-\big(\dot
c(s)\big)^*\Big)\beta(s)=$$
$$=\Big(\big\langle-,c_0(s)\big\rangle\sigma\dot
c_0(s)-\big\langle-,\dot c_0(s)\big\rangle\sigma
c_0(s)\Big)\beta(s)=\big\langle\beta(s),c_0(s)\big\rangle\sigma\dot
c_0(s)-\big\langle\beta(s),\dot c_0(s)\big\rangle\sigma c_0(s)$$
by Remark 2.2. It remains to observe that both $\beta(s)$ and $c_0(s)$
serve as solutions for the Cauchy problem
$\dot y(s)=\big\langle y(s),c_0(s)\big\rangle\sigma\dot
c_0(s)-\big\langle y(s),\dot c_0(s)\big\rangle\sigma c_0(s)$,
$y(a)=c_0(a)$, for an unknown $y:[a,b]\to V$
$_\blacksquare$

\medskip

{\bf2.5.~Lemma.} {\sl Let\/ $p\in\Bbb P_\Bbb CV\setminus\S V$, let\/
$c:[a,b]\to\Bbb P_\Bbb CV\setminus\S V$ be a piecewise smooth path, and
let\/ $c(a)=F_0p$ with\/ $F_0\in\SU V$. Then the isometry\/
$F:[a,b]\to\SU V$ that follows\/ $c$ and begins with\/ $F_0$ is a
horizontal lift of\/ $c$ with respect to the connection in the bundle\/
$\SU V\to\SU V/\Stab p$ whose horizontal spaces are the tangent
subspaces orthogonal to cosets with respect to the Killing form.}

\medskip

{\bf Proof.} For a nonisotropic $q$, the subspace
$\widehat{\T_q\Bbb P_\Bbb CV}\subset\su V$ coincides with
$(\T_1\Stab q)^\perp$. Indeed, let $l\in\T_1\Stab q$ and let
$t\in\T_q\Bbb P_\Bbb CV$, that is, $t=\langle-,q\rangle v$,
$\langle v,q\rangle=0$, $l\in\su V$, and $lq=0$. Then
$\tr(l\widehat t)=\langle lv,q\rangle=-\langle v,lq\rangle=0$. It
remains to compare the dimensions of $\widehat{\T_q\Bbb P_\Bbb CV}$ and
$(\T_1\Stab q)^\perp$.

In particular, $\dot F(s)\in\big(\T_1\Stab c(s)\big)^\perp\cdot F(s)$
for any $s\in[a,b]$. By Proposition 2.4, $F(s)=F_c(s)F_0$ and
$F(s)p=c(s)$. Hence, $\big(\Stab c(s)\big)F(s)=F(s)\Stab p$ and
$\dot F(s)\in F(s)\cdot(\T_1\Stab p)^\perp$
$_\blacksquare$

\medskip

We will describe isometries following geodesics. Due to Proposition
2.4, we should not worry much about the way of parameterizing
geodesics. So, we use the most convenient parameter. The {\it
bending\/} $B(s)$ {\it involving\/} $p_1,p_2$ is in fact the isometry
following the geodesic $\G{\wr}p_1,p_2{\wr}$ taken with a certain\break
parameterization proportional to the natural one. It is explicitly
defined in the proof of the following proposition.

\medskip

{\bf2.6.~Proposition.} {\sl Let\/
$p_1,p_2\in\Bbb P_\Bbb CV\setminus\S V$ be distinct and nonorthogonal.
All solutions\/ $x_1,x_2\in\Bbb P_\Bbb CV\setminus\S V$ of the
equation\/ $R(x_2)R(x_1)=R(p_2)R(p_1)$ can be described as follows.

There exists a one-parameter subgroup\/ $B:(\Bbb R,+)\to\SU V$ called\/
{\rm bending involving} $p_1,p_2$ such that\/ $B(s)$ commutes with\/
$R(p_2)R(p_1)$ and\/ $B(s)\G=\G$ for all\/ $s\in\Bbb R$, where\/
$\G:=\G{\wr}p_1,p_2{\wr}$ stands for the geodesic joining\/ $p_1,p_2$.
If\/ $p_1,p_2$ have the same sign, then\/ $B(1)p_1=p_2$. If\/ $p_1,p_2$
have different signs, the point\/ $B(1)p_1\in\G$ is orthogonal to\/
$p_2$.

For an arbitrary\/ $s\in\Bbb R$, we define\/ $q_j:=B(s)p_j$, $j=1,2$.
Then\/ $R(q_2)R(q_1)=R(p_2)R(p_1)$. If\/ $\G$ is noneuclidean, we
denote by\/ $q'_j\in\G$ the point orthogonal to\/ $q_j\in\G$, $j=1,2$.
Then\/ $R(q'_2)R(q'_1)=R(p_2)R(p_1)$.}

\medskip

{\bf Proof.} Denote by $p$ the polar point to $\L:=\L(p_1,p_2)$. Note
that the points $x_1,x_2$ subject to $R(x_2)R(x_1)=R(p_2)R(p_1)$ cannot
be equal or orthogonal by Remark 3.2. If $\L$ is noneuclidean, the
equality $R(q_2)R(q_1)=R(p_2)R(p_1)$ implies the equality
$R(q'_2)R(q'_1)=R(p_2)R(p_1)$. Indeed, using the orthogonal relations
(see Remark 3.2), we obtain
$R(p_2)R(p_1)=R(q_2)R(q_1)=R(q_2)R(p)R(p)R(q_1)=R(q'_2)R(q'_1)$.

\smallskip

$\bullet$ Suppose that $\G$ is hyperbolic. Choose representatives of
the vertices $v_1,v_2\in\G\cap\S V$ such that the Gram matrix of
$v_1,v_2,p$ equals
$\left[\smallmatrix0&\frac12&0\\\frac12&0&0\\0&0&1\endsmallmatrix
\right]$,
$p_1=v_1+\sigma_1v_2$, and $p_2=e^{-a}v_1+\sigma_2e^av_2$,
$0\ne a\in\Bbb R$, where $\sigma_j:=\sigma p_j$, $j=1,2$. The
positive/negative points in $\G$ have the form
$c_\pm(s):=e^{-as}v_1\pm e^{as}v_2$, $s\in\Bbb R$. It is easy to see
that $\big\langle c_\pm(s),c_\pm(s)\big\rangle=\pm1$ and
$\big\langle c_\pm(s),\dot c_\pm(s)\big\rangle=0$. So, the lifts are
normalized. By~Remark 2.2,
$\dot c_\pm(s)=\pm\big\langle-,c_\pm(s)\big\rangle\dot c_\pm(s)$. It
follows that
$\widehat{\dot c}(s):=\widehat{\dot c}_+(s)=\widehat{\dot
c}_-(s)=\left[\smallmatrix-a&0&0\\0&a&0\\0&0&0\endsmallmatrix\right]$.
Obviously,
$F(s):=F_{c_+}(s)=F_{c_-}(s)=\left[\smallmatrix
e^{-as}&0&0\\0&e^{as}&0\\0&0&1\endsmallmatrix\right]\in\SU(2,1)$
meets the conditions $\dot F(s)=\widehat{\dot c}(s)F(s)$ and $F(0)=1$.

We define $B(t):=F(t)$ and call the one-parametric subgroup
$B:(\Bbb R,+)\to\SU V$ the {\it bending involving\/} $p_1,p_2$. Note
that $B(1)p_1=p_2$ if $\sigma_1=\sigma_2$ and $B(1)p_1\in\L$ is
orthogonal to $p_2\in\L$ if $\sigma_1\ne\sigma_2$.

It is easy to verify that $B(s')c_\pm(s)=c_\pm(s+s')$ and
$R\big(c_\pm(s)\big)=\left[\smallmatrix0&\pm
e^{-2as}&0\\\pm e^{2as}&0&0\\0&0&-1\endsmallmatrix\right]$
for all $s,s'\in\Bbb R$. So,
$R\big(c_{\sigma_2}(s_2)\big)R\big(c_{\sigma_1}(s_1)\big)=
\left[\smallmatrix\sigma_1\sigma_2e^{-2a(s_2-s_1)}&0&0\\
0&\sigma_1\sigma_2e^{2a(s_2-s_1)}&0\\0&0&1\endsmallmatrix\right]$
and $R(q_2)R(q_1)=R(p_2)R(p_1)$.

For the converse, it suffices to observe that the points
$x_1,x_2\in\Bbb P_\Bbb CV\setminus\S V$ meeting the equation
$R(x_2)R(x_1)=R(p_2)R(p_1)$ have to belong to the geodesic $\G$ because
$\G$ is determined by the fixed points $v_1,v_2,p$ of $R(p_2)R(p_1)$.

\smallskip

$\bullet$ Suppose that $\G$ is spherical. Let $p'_1\in\G$ be such that
$p_1,p'_1$ are orthonormal. For every fixed $0\ne a\in\Bbb R$, the
points in $\G$ have the form $c(s):=p_1\cos as+p'_1\sin as$,
$s\in\Bbb R$. It is easy to see that this formula provides a normalized
lift of $c$. Choose $a\in(0,\frac{\pi}{2})$ and representatives of
$p_1,p'_1,p_2$ such that $p_1,p'_1$ are orthonormal and $c(1)=p_2$. By
Remark 2.2, $\dot c(s)=\big\langle-,c(s)\big\rangle\dot c(t)$. In the
orthonormal basis $p_1,p'_1,p$, we have
$\widehat{\dot
c}(s)=\left[\smallmatrix0&-a&0\\a&0&0\\0&0&0\endsmallmatrix\right]$.
So,
$F_c(s)=\left[\smallmatrix\cos as&-\sin as&0\\\sin as&\cos
as&0\\0&0&1\endsmallmatrix\right]\in\SU(2,1)$
meets the conditions $\dot F_c(s)=\widehat{\dot c}(s)F_c(s)$ and
$F_c(0)=1$.

We define $B(s):=F_c(s)$ and call the one-parametric subgroup
$B:(\Bbb R,+)\to\SU V$ the {\it bending involving\/} $p_1,p_2$.
Obviously, $B(1)p_1=p_2$.

It is easy to verify that $B(s')c(s)=c(s+s')$ and
$R\big(c(s)\big)=\left[\smallmatrix\cos2as&\sin2as&0\\
\sin2as&-\cos2as&0\\0&0&-1\endsmallmatrix\right]$.
This implies that
$R\big(c(s_2)\big)R\big(c(s_1)\big)=
\left[\smallmatrix\cos2a(s_2-s_1)&-\sin2a(s_2-s_1)&0\\
\sin2a(s_2-s_1)&\cos2a(s_2-s_1)&0\\0&0&1\endsmallmatrix\right]$
and $R(q_2)R(q_1)=R(p_2)R(p_1)$.

It remains to observe that the fixed points $p_1+ip'_1,p_1-ip'_1,p$ of
$R(p_2)R(p_1)$ determine\footnote{I am indebted to Carlos Henrique
Grossi Ferreira for pointing out this nice and easy fact.}
the geodesic $\G$ since
$\G=\big\{g\in\L\mid\ta(g,p+ip')=\ta(g,p-ip')\big\}$.

\smallskip

$\bullet$ Suppose that $\G$ is euclidean. Choose
$b\in\B V\cap\Bbb P_\Bbb Cp_1^\perp$ and representatives of
$b,p_1,p,p_2$ such that the Gram matrix of $b,p_1,p$ equals
$\left[\smallmatrix-1&0&1\\0&1&0\\1&0&0\endsmallmatrix\right]$ and
$p_2=p_1+ap$ for some $a>0$. Every point in $\G$ but $p$ has the form
$c_0(s):=p_1+asp$, $s\in\Bbb R$. It is easy to see that another choice
of $b$ and the representatives does not change the parameterization
$c(s)$ of $\G$ (however, it may change $a$). Since
$\big\langle c_0(s),c_0(s)\big\rangle=1$ and
$\big\langle c_0(s),\dot c_0(s)\big\rangle=0$, the lift $c_0$ is
normalized. By Remark 2.2,
$\dot c(s)=\big\langle-,c_0(s)\big\rangle\dot
c_0(s)=\big\langle-,c_0(s)\big\rangle ap$.
It~follows that
$\widehat{\dot c}_0(s)=\left[\smallmatrix0&0&0\\-a&0&0\\0&a&0
\endsmallmatrix\right]$
in the basis $b,p_1,p$. Therefore,
$F_c(s)=\left[\smallmatrix1&0&0\\-as&1&0\\-a^2s^2/2&as&1\endsmallmatrix
\right]\in\SU(2,1)$
meets the conditions $\dot F_c(s)=\widehat{\dot c}(s)F_c(s)$ and
$F_c(0)=1$.

We define $B(s):=F_c(s)$ and call the one-parametric subgroup
$B:(\Bbb R,+)\to\SU(2,1)$ the {\it bending involving\/} $p_1,p_2$.
Obviously, $B(1)p_1=p_2$.

A straightforward verification shows that $B(s')c_0(s)=c_0(s+s')$ and
$R\big(c_0(s)\big)=\left[\smallmatrix-1&0&0\\2as&1&0\\2a^2s^2&2as&-1
\endsmallmatrix\right]$.
Hence,
$R\big(c_0(s_2)\big)R\big(c_0(s_1)\big)=\left[\smallmatrix1&0&0\\-
2a(s_2-s_1)&1&0\\-2a^2(s_2-s_1)^2&2a(s_2-s_1)&1\endsmallmatrix\right]$.
So, $R(q_2)R(q_1)=R(p_2)R(p_1)$ and these are unique solutions
$x_1,x_2\in\G$ of the equation $R(x_2)R(x_1)=R(p_2)R(p_1)$.

As $p$ is the only fixed point of $R(p_2)R(p_1)$, the points $x_1,x_2$
subject to the equation $R(x_2)R(x_1)=R(p_2)R(p_1)$ should belong to
$\L$. The geodesic $\G$ is stable under the action of $R(p_2)R(p_1)$.
So is the geodesic $\G{\wr}x_1,x_2{\wr}\subset\L$. Writing down the
$x_j$'s in the form $x_j=p_1+c_jp$, $c_j\in\Bbb C$, we conclude that
$p_1+(2a+c_j)p\in\G{\wr}p_1+c_1p,p_1+c_2p{\wr}$ and derive from this
that $c_1,c_2\in\Bbb R$
$_\blacksquare$

\medskip

It is not difficult to show that the slice identification introduced in
[AGG] and the meridional transport mentioned in [AGr] are induced by
suitable bendings. Also, one can derive the formulae for parallel
transport along geodesics, given in [AGr], from the above description
of bendings.

\bigskip

\centerline{\bf3.~Bendings and reflections}

\medskip

In this section, we study how an arbitrary relation between reflections
can be modified with the help of short relations ($n\le4$). Also, we
show that almost all isometries in $\SU V$ are products of $3$
reflections.

\smallskip

Let $p\in\Bbb P_\Bbb CV\setminus\S V$. The {\it reflection\/}
$R(p):V\to V$ in $p$ is defined as
$$R(p):x\mapsto2\textstyle\frac{\langle x,p\rangle}{\langle
p,p\rangle}p-x.$$
It is easy to see that $R(p)\in\SU V$ and $R(p)R(p)=1$.

\smallskip

Given pairwise orthogonal $p_1,p_2,p_3\in\Bbb P_\Bbb CV\setminus\S V$,
it is immediate that the $R(p_i)$'s commute and $R(p_3)R(p_2)R(p_1)=1$.

\smallskip

Denote
$$\ta(p_1,p_2):=\textstyle\frac{g_{12}g_{21}}{g_{11}g_{22}},\ \
\alpha(p_1,p_2,p_3):=\Im\frac{g_{12}g_{23}g_{31}}{g_{11}g_{22}g_{33}},\
\ \beta(p_1,p_2,p_3):=\frac{\det[g_{jl}]}{g_{11}g_{22}g_{33}},\ \
\tau(p_1,p_2,p_3):=\Re\frac{g_{13}g_{22}}{g_{12}g_{23}},$$
where $[g_{jk}]$ stands for the Gram matrix of nonisotropic
$p_1,p_2,p_3\in V$, i.e., $g_{jk}:=\langle p_j,p_k\rangle$.

\medskip

{\bf3.1.~Lemma {\rm(compare with [Pra])}.} {\sl Let\/ $[g_{jk}]$ be the
Gram matrix of\/ $p_1,\dots,p_n\in\Bbb P_\Bbb CV\setminus\S V$. Then
$$\tr\big(R(p_n)\dots R(p_1)\big)=(-1)^n\Big(3-2n+\sum\limits_{1\le
i_1<\dots<i_t\le n\atop2\le t\le
n}(-2)^t\frac{g_{i_1i_2}g_{i_2i_3}\dots
g_{i_{t-1}i_t}g_{i_ti_1}}{g_{i_1i_1}g_{i_2i_2}\dots
g_{i_{t-1}i_{t-1}}g_{i_ti_t}}\Big).$$
In particular,
$$\tr\big(R(p_1)\big)=-1,\qquad\tr\big(R(p_2)R(p_1)\big)=
4\ta(p_1,p_2)-1,$$
$$\tr\big(R(p_3)R(p_2)R(p_1)\big)=8i\alpha(p_1,p_2,p_3)+
4\beta(p_1,p_2,p_3)-1.$$}
\indent
{\bf Proof.} The facts that the trace of $x\mapsto\langle x,q\rangle p$
equals $\langle p,q\rangle$ and that $\tr1=3$ imply the general
formula. The rest follows from
$$\det[g_{jl}]=g_{11}g_{22}g_{33}+2\Re(g_{12}g_{23}g_{31})-
g_{13}g_{31}g_{22}-g_{12}g_{21}g_{33}-g_{23}g_{32}g_{11}=$$
$$=g_{11}g_{22}g_{33}\big(1+2\textstyle\frac{\Re(g_{12}g_{23}g_{31})}
{g_{11}g_{22}g_{33}}-\ta(p_1,p_2)-\ta(p_2,p_3)-\ta(p_3,p_1)\big),$$
$$\tr\big(R(p_3)R(p_2)R(p_1)\big)=3-4\big(\ta(p_1,p_2)+
\ta(p_2,p_3)+\ta(p_3,p_1)\big)+8\textstyle\frac{g_{12}g_{23}g_{31}}
{g_{11}g_{22}g_{33}}\ _\blacksquare$$

The following well-known facts are straightforward or easily follow
from Lemma 3.1 :

\medskip

{\bf3.2.~Remark.} {\sl Let\/ $\delta\in\Bbb C$ be such that\/
$\delta^3=1$ and let\/ $p_1,p_2,p_3\in\Bbb P_\Bbb CV\setminus\S V$.

\smallskip

$\bullet$ If\/ $F^2=\delta$ for some\/ $F\in\SU V$, then either\/
$F=\delta^2$ or\/ $F=\delta^2R(p)$ for a suitable\/
$p\in\Bbb P_\Bbb CV\setminus\S V$.

$\bullet$ If\/ $R(p_2)R(p_1)=\delta$, then\/ $\delta=1$ and\/
$p_1=p_2$.

$\bullet$ If\/ $R(p_3)R(p_2)R(p_1)=\delta$, then\/ $\delta=1$ and\/
$p_1,p_2,p_3$ are pairwise orthogonal
$_\blacksquare$}

\medskip

Suppose that $R(q_1)R(p_1)=\delta R(q_2)R(p_2)$ with $\delta^3=1$. Then
$\delta=1$. Indeed, $4\ta(p_1,q_1)-1=\delta\big(4\ta(p_2,q_2)-1\big)$
by Lemma 3.1. If $\delta\ne1$, then $4\ta(p_j,q_j)=1$ and the
projective line $\L(p_j,q_j)$ is spherical, $j=1,2$. Hence, the polar
point $b_j$ to $\L(p_j,q_j)$ is negative. Being $b_j$ the only negative
fixed point of $R(q_j)R(p_j)$, we obtain $b_1=b_2$. The other two fixed
points of $R(q_1)R(p_1)$ and of $R(q_2)R(p_2)$ also coincide. For the
reason used in the proof of Proposition 2.6,
$\G{\wr}p_1,q_1{\wr}=\G{\wr}p_2,q_2{\wr}$. Applying a suitable bending
involving $p_2,q_2$, we can assume that $p_2=p_1$. A contradiction.

\smallskip

So, from Remark 3.2 and Proposition 2.6, we obtain the following
corollary.

\medskip

{\bf3.3.~Corollary.} {\sl All relations with\/ $n\le4$ follow from
cancellations, orthogonal relations, and bending relations
$_\blacksquare$}

\medskip

{\bf3.4.~Remark.} {\sl Given\/
$p_1,p_2,p_3\in\Bbb P_\Bbb CV\setminus\S V$ such that\/ $p_1,p_2$ are
distinct nonorthogonal and\/ $\L:=\L(p_1,p_2)$ is not spherical, the
projective line\/ $\L(p_2,p_3)$ becomes hyperbolic after a suitable
bending involving\/ $p_1,p_2$ except in the following cases\/{\rm:}

\smallskip

$\bullet$ $\L$ is euclidean and\/ $p_3\in\L$,

$\bullet$ $\L$ is hyperbolic and\/ $p_3$ is the polar point to\/ $\L$.}

\medskip

{\bf Proof.} Denote by $p$ the polar point to $\L$. We can assume that
$\langle p_2,p_2\rangle=\langle p_3,p_3\rangle=1$.

Suppose that $\L$ is euclidean and that $p_3\notin\L$, i.e.,
$\langle p,p_3\rangle\ne0$. Every point in
$\G{\wr}p_1,p_2{\wr}\setminus\{p\}$ has the form $p_2+sp$,
$s\in\Bbb R$, after a suitable choice of a representative $p\in V$.
Hence,
$\ta(p_2+sp,p_3)=\big|\langle p_2,p_3\rangle+s\langle
p,p_3\rangle\big|^2\to\infty$
while $s\to\infty$ because $\langle p,p_3\rangle\ne0$.

Suppose that $\L$ is hyperbolic. Every positive point in
$\G{\wr}p_1,p_2{\wr}$ has the form $p(s)=e^{-s}v_1+e^sv_2$,
$s\in\Bbb R$, where $v_1,v_2$ stand for representatives of the vertices
of $\G{\wr}p_1,p_2{\wr}$ such that $\langle v_1,v_2\rangle=\frac12$.
Therefore,
$\ta\big(p(t),p_3\big)=\big|e^{-s}\langle v_1,p_3\rangle+e^s\langle
v_2,p_3\rangle\big|^2\to\infty$
while $s\to\pm\infty$ unless
$\langle v_1,p_3\rangle=\langle v_2,p_3\rangle=0$
$_\blacksquare$

\medskip

{\bf3.5.~Corollary.} {\sl Let\/ $R(p_n)\dots R(p_1)=\delta$, where\/
$p_1,\dots,p_n\in\Bbb P_\Bbb CV\setminus\S V$ and\/ $\delta^3=1$. Then,
using cancellations, orthogonal relations, and bendings, we can either
diminish $n$ or reach a situation where at most one of the\/ $p_j$'s is
positive unless every projective line\/ $\L(p_{j-1},p_j)$ is
spherical\/ {\rm(}the indices are modulo\/ $n${\rm)}.}

\medskip

{\bf Proof.} First, we show how to gain a negative $p_j$. If some
$\L(p_{j-1},p_j)$ is hyperbolic, this is easy since, using orthogonal
relations, we can simultaneously alter the signs of $p_{j-1},p_j$, if
necessary. So, we assume $\L(p_1,p_2)$ euclidean. By Remark 3.4, we
assume that $p_3\in\L(p_1,p_2)$, that is, $\L(p_1,p_2)=\L(p_2,p_3)$.
In~this way, we arrive at a situation where the $p_j$'s belong to the
same euclidean projective line $\L$. It~remains to apply elementary
euclidean plane arguments. (See the end of [AGr, Subsection 3.1] for
the description of geodesics in an euclidean projective line.) If the
geodesics $\G{\wr}p_1,p_2{\wr}$ and $\G{\wr}p_3,p_4{\wr}$ intersect,
i.e.,\break

\vskip-8pt

\noindent
$\vcenter{\hbox{\epsfbox{Picture1.eps}}}$\hskip10pt
$\vcenter{\hbox{\epsfbox{Picture2.eps}}}$

\leftskip127pt

\vskip-47pt

\noindent
if they have a common point $q$ different from the polar point to $\L$,
then we obtain $p_3=q$ after a suitable bending involving $p_3,p_4$.
Hence, $p_3\in\G{\wr}p_1,p_2{\wr}$ and some bending involving $p_2,p_3$
provides $p_2=p_1$, i.e., a cancellation. If
$p_5\in\G{\wr}p_3,p_4{\wr}$, then we can get $p_3=p_4$ with a bending
involving\break

\vskip-5pt

\leftskip0pt

\noindent
\hskip396pt$\vcenter{\hbox{\epsfbox{Picture3.eps}}}$

\rightskip62pt

\vskip-48pt

\noindent
$p_4,p_5$. If $p_5\notin\G{\wr}p_3,p_4{\wr}$, a `small' bending
involving $p_4,p_5$ makes $\G{\wr}p_1,p_2{\wr}$ and
$\G{\wr}p_3,p_4{\wr}$ intersect.

Without loss of generality, we assume $p_1,\dots,p_{j-1}$ negative and
$p_j,p_n$ positive, $1<j<\nomathbreak n$. If $p_{j+1}$ is the polar
point to $\L(p_{j-1},p_j)$, we can diminish $n$ by using an
orthogonal\break

\vskip-12pt

\rightskip0pt

\noindent
relation. Otherwise, by Remark 3.4, $\L(p_j,p_{j+1})$ becomes
hyperbolic after a bending involving $p_{j-1},p_j$. It remains to alter
the signs of $p_j,p_{j+1}$, if necessary
$_\blacksquare$

\medskip

{\bf3.6.~Definition.} Let $d=2$. A triple of points
$p_1,p_2,p_3\in\Bbb P_\Bbb CV\setminus\S V$ is said to be {\it
regular\/} if $p_1,p_2,p_3$ are not in a same geodesic, $p_2$ is not
orthogonal to $p_j$ for $j=1,3$, and at most one of $p_1,p_2,p_3$ is
positive. A regular triple $p_1,p_2,p_3$ is {\it strongly regular\/}
when $p_1,p_2,p_3$ are not in a same projective line and $p_1,p_2,p_3$
are all negative if they belong to a same real plane. An isometry
$F\in\SU V$ is called {\it regular\/} if
$\dim_\Bbb C\{v\in V\mid Fv=cv\}\le1$ for any $c\in\Bbb C$.

\medskip

{\bf3.7.~Lemma.} {\sl Let\/ $p_1,p_2,p_3$ be a regular triple. Then the
isometry\/ $F:=R(p_3)R(p_2)R(p_1)$ is regular.}

\medskip

{\bf Proof.} Suppose that $\dim_\Bbb C\{v\in V\mid Fv=cv\}\ge2$. Then
$|c|=1$. Note that $c\ne-1$. Otherwise, $F=R(p)$ is a reflection and
$p,p_1,p_2,p_3$ are in a same geodesic by Proposition 2.6,
contradicting the assumption that the triple is regular. We can pick
$0\ne v\in V$ such that $\langle v,p_1\rangle=0$ and $Fv=cv$. This
means that $R(p_1)v=-v$ and, hence, $R(p_2)v+cR(p_3)v=0$. In other
words, assuming that $\sigma_j:=\langle p_j,p_j\rangle\in\{-1,1\}$, we
obtain
$2\sigma_2\langle v,p_2\rangle p_2+2\sigma_3c\langle v,p_3\rangle
p_3=(1+c)v$.
Since $1+c\ne0$ and $p_2,p_3$ are $\Bbb C$-linearly independent, we
conclude that $v\in\Bbb Cp_2+\Bbb Cp_3$. If $v$ is proportional to one
of $p_2,p_3$, then $\langle p_2,p_3\rangle=0$, which is impossible for
a regular triple. Hence, we can assume that $v=zp_2+p_3$, for some
$0\ne z\in\Bbb C$. By a straightforward calculation, we obtain
$\langle p_2,p_3\rangle=\frac{\sigma_2(\overline c-1)\overline
z}2=\frac{\sigma_3(1-c)}{2cz}$.
Therefore, $|z|=1$ and $\sigma_2=\sigma_3$, implying
$\sigma_2=\sigma_3=-1$ and $\big|\langle p_2,p_3\rangle\big|<1$. It
follows that $\ta(p_2,p_3)<1$, a~contradiction
$_\blacksquare$

\medskip

{\bf3.8.~Proposition.} {\sl An isometry \/ $F\in\SU V$ admits the
form\/ $F=R(p_3)R(p_2)R(p_1)$, where\/ $p_1,p_2,p_3$ is a regular
triple, iff\/ $F$ is regular and\/ $\tr F\ne-1$.}

\medskip

{\bf Proof.} If $\tr\big(R(p_3)R(p_2)R(p_1)\big)=-1$, then
$\alpha(p_1,p_2,p_3)=\beta(p_1,p_2,p_3)=0$ by Lemma 3.1. Hence,
$p_1,p_2,p_3$ are in a same geodesic.

Conversely, given a regular $F\in\SU V$ such that $\tr F\ne-1$, write
$\tr F=8i\alpha+4\beta-1$ and $\beta=\mp|\beta|$, where
$\alpha,\beta\in\Bbb R$. For sufficiently big $g_{12},g_{23}>1$, we
find $t\in\Bbb R$ satisfying the equation
$$g_{12}^2g_{23}^2(t-1)^2=(g_{12}^2\pm1)(g_{23}^2\pm1)-\beta-
\textstyle\frac{\alpha^2}{g_{12}^2g_{23}^2}.$$
We have
$$\pm\det G_\pm=1+2g_{12}^2g_{23}^2t\pm g_{12}^2\pm
g_{23}^2-g_{12}^2g_{23}^2t^2-\textstyle\frac{\alpha^2}
{g_{12}^2g_{23}^2}=(g_{12}^2\pm1)(g_{23}^2\pm1)-g_{12}^2g_{23}^2(t-1)^2
-\frac{\alpha^2}{g_{12}^2g_{23}^2}=\beta=\mp|\beta|,$$
where
$$G_\pm:=\left[\smallmatrix-1&g_{12}&\pm g_{12}g_{23}t\mp\frac{\alpha
i}{g_{12}g_{23}}\\g_{12}&\pm1&g_{23}\\\pm g_{12}g_{23}t\pm\frac{\alpha
i}{g_{12}g_{23}}&g_{23}&-1\endsmallmatrix\right].$$
By Sylvester's Criterion, there exist $p_1,p_2,p_3\in V$ with the Gram
matrix $G_\pm$ because $g_{12}>1$ and $\det G_\pm\le0$. Note that
$p_1,p_2,p_3$ are not in the same geodesic because otherwise
$\alpha=\beta=0$ and $\tr F=-1$. It remains to apply Lemma 3.1, [Gol2,
Theorem 6.2.4], and Lemma 3.7
$_\blacksquare$

\medskip

In the sequel, we will need the following well-known fact which can be
obtained by straightforward calculations:

\medskip

{\bf3.9.~Remark.} {\sl Let\/ $F\in\SU V$ be a regular isometry, let\/
$C_0(F)\le\U V$ stand for the centralizer of\/ $F$ in\/ $\U V$, and
let\/ $C(F)\le\PU V$ denote the image of\/ $C_0(F)$. Then\/ $C(F)$ is
connected, commutative, and\/~$2$-dimensional.}
$_\blacksquare$

\bigskip

\centerline{\bf4.~Composition of bendings}

\medskip

Let $p_1,p_2,p_3$ be a strongly regular triple (see Definition 3.6). In
this section, we study compositions of two types of bendings: those
involving $p_1,p_2$ and those involving $p_2,p_3$.

\smallskip

These bendings preserve the regular isometry $F:=R(p_3)R(p_2)R(p_1)$
and the numbers $\sigma_j:=\sigma p_j$, $j=1,2,3$. So, we fix the
$\sigma_j$'s and $\alpha,\beta\in\Bbb R$ such that
$\tr F=8i\alpha+4\beta-1$, $\alpha:=\alpha(p_1,p_2,p_3)$,
and~$\beta:=\beta(p_1,p_2,p_3)$ (see Lemma 3.1). By Definition 3.6,
Lemma 3.1, Proposition 3.8, and Sylvester's Criterion, the only
restrictions on the numbers $\sigma_1,\sigma_2,\sigma_3,\alpha,\beta$
are

\smallskip

$\bullet$ $\sigma_1\sigma_2\sigma_3\beta<0$, at most one of the
$\sigma_j$'s is positive, and
$\alpha=0\Longrightarrow\sigma_1=\sigma_2=\sigma_3=-1$.

\smallskip

We will frequently choose representatives of $p_1,p_2,p_2$ such that
$g_{jj}=\sigma_j$ and $g_{12},g_{23}>0$, where $[g_{jk}]$ stands for
the Gram matrix of $p_1,p_2,p_2$. Such a Gram matrix is called {\it
standard.} Note that the standard Gram matrix is uniquely determined by
a strongly regular triple $p_1,p_2,p_2\in\Bbb P_\Bbb CV$.

When $\alpha=0$, the standard Gram matrix is real. In this case, the
triple lies in (and spans) a real plane. Such a triple is said to be
{\it real.}

\medskip

{\bf4.1.~Geometrical configurations of strongly regular triples.}
First, we will show that, geometrically, i.e., up to the action of
$\PU V$, the strongly regular triples $p_1,p_2,p_3$ with fixed
$\sigma_1,\sigma_2,\sigma_3$ and the conjugacy class of
$F:=R(p_3)R(p_2)R(p_1)$ form a surface $S$. Then we describe the lines
on $S$ that correspond to the bendings involving $p_1,p_2$/$p_2,p_3$.

\medskip

{\bf4.1.1.~Lemma.} {\sl Geometrically, a strongly regular triple\/
$p_1,p_2,p_3$ can be described as a point in the surface\/
$S\subset\Bbb R^3(t,t_1,t_2)$ given by the equation
$$(t_1-1)(t_2-1)=t_1t_2(t-1)^2+\textstyle\frac{\alpha^2}{t_1t_2}+
\beta\leqno{\bold{(4.1.2)}}$$
and by the inequalities
$$\sigma_1\sigma_2t_1>0,\qquad\sigma_1\sigma_2t_1>\sigma_1\sigma_2,
\qquad\sigma_2\sigma_3t_2>0,\qquad\sigma_2\sigma_3t_2>
\sigma_2\sigma_3,\leqno{\bold{(4.1.3)}}$$
where\/ $t_1:=\ta(p_1,p_2)$, $t_2:=\ta(p_2,p_3)$, and\/
$t:=\tau(p_1,p_2,p_3):=\Re\frac{g_{13}g_{22}}{g_{12}g_{23}}$.}

\medskip

{\bf Proof.} By Sylvester's Criterion, a strongly regular triple
$p_1,p_2,p_3$ is geometrically given by its Gram matrix $[g_{jk}]$. We
assume it to be standard. Let $t_1:=\ta(p_1,p_2)$, $t_2:=\ta(p_2,p_3)$,
and $t:=\tau(p_1,p_2,p_3):=\Re\frac{g_{13}g_{22}}{g_{12}g_{23}}$. It
follows from $\alpha(p_1,p_2,p_3)=\alpha$ that
$g_{31}=\frac{g_{12}g_{23}}{g_{22}}t+
\frac{g_{11}g_{22}g_{33}}{g_{12}g_{23}}\alpha i$.
Therefore,
$$\beta:=\textstyle\frac{\det[g_{jk}]}{g_{11}g_{22}g_{33}}=1+2t_1t_2
-t_1-t_2-\frac{g_{12}^2g_{23}^2}{g_{11}g_{22}^2g_{33}}t^2-
\frac{g_{11}g_{22}^2g_{33}}{g_{12}^2g_{23}^2}\alpha^2=
1+2t_1t_2-t_1-t_2-t_1t_2t^2-\frac{\alpha^2}{t_1t_2}$$
and we arrive at (4.1.2). The inequalities
$\sigma_1\sigma_2t_1>\sigma_1\sigma_2$ and
$\sigma_2\sigma_3t_2>\sigma_2\sigma_3$ follow from the fact that the
projective lines $\L(p_1,p_2)$ and $\L(p_2,p_3)$ are hyperbolic.

Conversely, let $(t,t_1,t_2)\in S$. We put $g_{jj}:=\sigma_j$ for
$j=1,2,3$, $g_{12}:=g_{21}:=\sqrt{\sigma_1\sigma_2t_1}>0$,
$g_{23}:=g_{32}:=\sqrt{\sigma_2\sigma_3t_2}>0$,
$g_{31}:=\frac{g_{12}g_{23}}{g_{22}}t+
\frac{g_{11}g_{22}g_{33}}{g_{12}g_{23}}\alpha i$,
and $g_{13}:=\overline g_{31}$. The equation (4.1.2) implies that
$\det[g_{jk}]=\sigma_1\sigma_2\sigma_3\beta$. By Sylvester's Criterion,
there are points $p_1,p_2,p_3\in V$ with the Gram matrix $[g_{jk}]$
because at most one of the $\sigma_j$'s is positive,
$\sigma_1\sigma_2t_1>\sigma_1\sigma_2$, and
$\sigma_1\sigma_2\sigma_3\beta<0$
$_\blacksquare$

\medskip

Given $r_1,r_2\in\Bbb R$ such that $\sigma_1\sigma_2r_1>0$,
$\sigma_1\sigma_2r_1>\sigma_1\sigma_2$, $\sigma_2\sigma_3r_2>0$, and
$\sigma_2\sigma_3r_2>\sigma_2\sigma_3$, denote by $V_{r_1}$, $H_{r_2}$,
and $C$ the lines given in $S$ respectively by the equations $t_1=r_1$,
$t_1=r_2$, and $t=1$. We call $V_{r_1}$/$H_{r_2}$ a {\it
vertical\/}/{\it horizontal\/} line. Obviously, the intersection of a
vertical line and a horizontal one is either empty or a point in $C$ or
consists of two points outside $C$.

\vskip19pt

\noindent
\hskip270pt$\vcenter{\hbox{\epsfbox{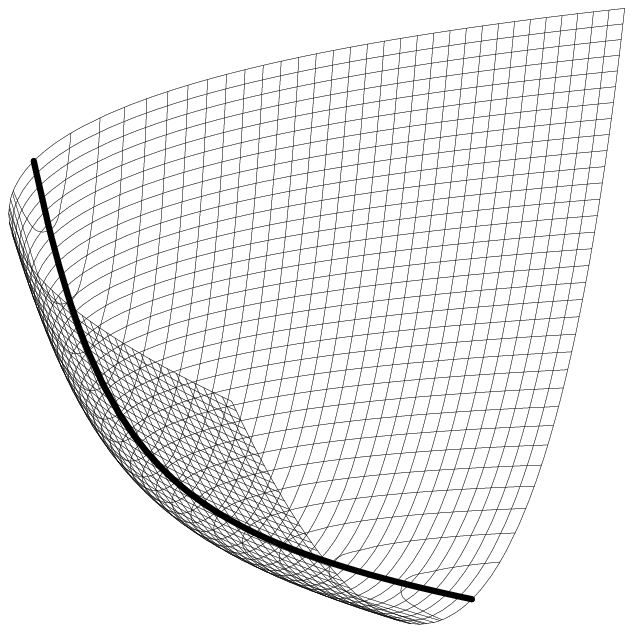}}}$

\rightskip195pt

\vskip-198pt

\medskip

{\bf4.1.4.~Lemma.} {\sl Every vertical\/{\rm/}horizontal line is
non\-empty, smooth, and connected. It intersects\/ $C$ in exactly\/ $1$
point and corresponds to the bendings involving the points\/
$p_1,p_2${\rm/}$p_2,p_3$. The surface\/ $S$ is a topological plane
fibred by vertical\/{\rm/}horizontal lines. If a vertical and a
horizontal lines intersect outside\/~$C$, this intersection is
transversal.}

\medskip

{\bf Proof.} Let $Q\subset\Bbb R^2(t_1,t_2)$ denote the quadrant given
by the inequalities (4.1.3). The projection of the surface $S$ into the
plane $\Bbb R^2(t_1,t_2)$ coincides with the region $R\subset Q$ given
in $Q$ by the inequality $f(t_1,t_2)\ge0$, where
$f(t_1,t_2):=\sigma_1\sigma_3\big((t_1-1)(t_2-1)-
\frac{\alpha^2}{t_1t_2}-\beta\big)$.
The map $S\to R$ is a double covering ramified along $C$. The fact that
vertical/horizontal lines of $S$ are smooth follows
straightforwardly.\footnote{For example, if the derivatives of
$r_1t_2(t-1)^2+\frac{\alpha^2}{r_1t_2}+\beta-(r_1-1)(t_2-1)$ with
respect to $t$ and $t_2$ equal respectively
$r_2(t-1)^2-\frac{\alpha^2}{r_1t_2^2}-(r_1-1)$ and $2r_1t_2(t-1)$. If
both vanish, then $t=1$ and $\frac{\alpha^2}{t_2^2}+r_1(r_1-1)=0$. On
the other hand,
$r_1(r_1-1)=
\sigma_1\sigma_2r_1(\sigma_1\sigma_2r_2-\sigma_1\sigma_2)>0$.}
Such lines project into the intersections with $R$ of the
vertical/horizontal lines of the plane $\Bbb R^2(t_1,t_2)$.

\rightskip0pt

Let $\sigma_2\sigma_3r_2>0$ and $\sigma_2\sigma_3r_2>\sigma_2\sigma_3$.
The function
$$f(\sigma_1\sigma_2x,r_2)=
(x-\sigma_1\sigma_2)(\sigma_2\sigma_3r_2-\sigma_2\sigma_3)-
\textstyle\frac{\alpha^2}{\sigma_2\sigma_3r_2x}+\sigma_2|\beta|$$
is increasing in $x>0$ and
$\lim\limits_{x\to+\infty}f(\sigma_1\sigma_2x,r_2)=+\infty$. Also, it
takes negative values for $x>\max(\sigma_1\sigma_2,0)$. Indeed, if
$\sigma_1\sigma_2=1$, then
$\lim\limits_{x\to+1}f(\sigma_1\sigma_2x,r_2)=
-\frac{\alpha^2}{\sigma_2\sigma_3r_2}-|\beta|$.
If $\sigma_1\sigma_2=-1$, then $\alpha\ne0$ and
$\lim\limits_{x\to+0}f(\sigma_1\sigma_2x,r_2)=-\infty$. We~conclude
that every horizontal line of $Q$ intersects $C$ in exactly $1$ point.
Hence, every horizontal line of $S$ is connected and intersects $C$ in
exactly $1$ point. By symmetry, the~same is true for vertical lines.

It remains to prove the claim concerning the bendings. By symmetry, we
deal only with vertical lines. Obviously, when we bend $p_1,p_2$, we
keep $t_1$ constant, $t_1=r_1$. Denote by $v_1,v_2$ the vertices of the
geodesic $\G{\wr}p_1,p_2{\wr}$ and choose representatives of
$p_1,p_2,p_3,v_1,v_2$ such that $\langle p_3,p_3\rangle=\sigma_3$,
$\langle v_1,v_2\rangle=\frac12$, $p_1=v_1+\sigma_1v_2$, and
$p_2=e^{-a}v_1+\sigma_2e^av_2$ with $a>0$ and $a\ne1$. The formulae
$$p_1(s):=e^{-as}v_1+\sigma_1e^{as}v_2,\qquad
p_2(s):=e^{-a(s+1)}v_1+\sigma_2e^{a(s+1)}v_2,\qquad s\in\Bbb R,$$
describe the bendings involving $p_1,p_2$. Denote
$z_j:=\langle v_j,p_3\rangle$ and $c_j:=|z_j|$ for $j=1,2$. Note that
$z_j\ne0$. Indeed, suppose that $z_2=0$. Then $\sigma_3=1$, hence,
$\sigma_1=\sigma_2=-1$ and $\alpha\ne0$. On the other hand,
$\alpha=\sigma_1\sigma_2\sigma_3\Im\big(\frac12(\sigma_1e^{-a}+
\sigma_2e^a)e^{-a(2s+1)}z_1\overline z_1\big)=0$,
a contradiction.

From
$$t(s)=\Re\textstyle\frac{2(e^{-as}z_1+\sigma_1e^{as}z_2)\sigma_2}
{(\sigma_1e^{-a}+\sigma_2e^a)(e^{-a(s+1)}z_1+\sigma_2e^{a(s+1)}z_2)}=
\frac{2(\sigma_2e^{-a(2s+1)}c_1^2+(\sigma_1\sigma_2e^{-a}+e^a)
\Re(z_1\overline z_2)+\sigma_1e^{a(2s+1)}c_2^2)}{(\sigma_1e^{-a}+
\sigma_2e^a)|e^{-a(s+1)}z_1+\sigma_2e^{a(s+1)}z_2|^2},$$
we obtain
$$\lim_{s\to-\infty}t(s)+\lim_{s\to+\infty}t(s)=
\textstyle\frac{2\sigma_2}{(\sigma_1e^{-a}+\sigma_2e^a)e^{-a}}+
\frac{2\sigma_1}{(\sigma_1e^{-a}+\sigma_2e^a)e^a}=2.$$
So, $\lim_{s\to-\infty}t(s)-1$ and $\lim_{s\to+\infty}t(s)-1$ have
opposite signs. On the other hand, the function
$$\sigma_2\sigma_3t_2(s)=|e^{-a(s+1)}z_1+\sigma_2e^{a(s+1)}z_2|^2=
e^{-2a(s+1)}c_1^2+e^{2a(s+1)}c_2^2+2\sigma_2\Re(z_1\overline z_2)$$
in $s$ possesses exactly one minimum and
$\lim\limits_{s\to\pm\infty}\sigma_2\sigma_3t_2(s)=+\infty$. Taking
into account the above description of a vertical line and its
projection into $Q$, this completes the proof
$_\blacksquare$

\medskip

{\bf4.1.5.~Corollary.} {\sl The bending involving\/
$p_1,p_2${\rm/}$p_2,p_3$ keeps any strongly regular triple\/
$p_1,p_2,p_3$ strongly regular. Geometrically, all strongly regular
triples with fixed\/ $\alpha$, $\beta$, and\/ $\sigma_j$'s are
connected by means of such bendings}
$_\blacksquare$

\medskip

{\bf4.2.~Holonomy of strongly regular triples.} It can easily happen
that, after a few bendings involving $p_1,p_2$/$p_2,p_3$, we obtain
from a strongly regular triple $p_1,p_2,p_3$ a different triple
$p'_1,p'_2,p'_3$ which is geometrically the same. In order to study
this phenomena, we introduce some formal settings.

\smallskip

\vskip7pt

{\unitlength=1bp$$\latexpic(0,0)(224,13)\put(0,24){$E_F$}
\put(15,27){\vector(1,0){12}}\put(28,24){$E$}
\put(37,27){\vector(1,0){10}}\put(49,24){$\SU V$}
\put(12,21){\vector(1,-1){15}}\put(32,22){\vector(0,-1){12}}
\put(34,16){$\pi$}\put(28,0){$S$}\endlatexpic$$}

\leftskip79pt

\vskip-48pt

Denote by $E$ the manifold of all strongly regular triples
$p_1,p_2,p_3\in\Bbb P_\Bbb CV$ (the~conjugacy class of
$R(p_3)R(p_2)R(p_1)\in\SU V$ and the $\sigma_j$'s are fixed). To every
such triple, we can associate its standard Gram matrix; this matrix
determines a triple of representatives $p_1,p_2,p_3\in V$, unique up to
a unitary factor. We obtain the map\break

\vskip-12pt

\leftskip0pt

\noindent
$\pi:E\to S$. By Sylvester's Criterion, $\pi$ is a principal
$\PU V$-bundle.

To every triple $p:=(p_1,p_2,p_3)\in E$, we associate the isometry
$F_p:=R(p_3)R(p_2)R(p_1)\in\SU V$. The left action of $\PU V$ on $E$
induces the conjugation at the level of $F_p$. Indeed, let $g\in\U V$
be a representative of $g\in\PU V$. Then
$$F_{gp}=R(gp_3)R(gp_2)R(gp_1)=gR(p_3)R(p_2)R(p_1)g^{-1}=
gF_pg^{-1}.$$
In particular, the triples $p$ with fixed $F_p=F$ form a principal
$C(F)$-bundle $E_F\to S$.

By Lemma 4.1.4, every vertical/horizontal line in $S$ lifts to the
trajectory of the corresponding bending. This provides a connection
over $S\setminus C$ because at any point $q\in S\setminus C$ the
tangent vectors to the vertical and horizontal lines passing through
$q$ generate the tangent space $\T_qS$. (However, at a point in $C$,
the vertical and horizontal lines are tangent.) Denote by $H(q_0)$ the
corresponding (connected) holonomy group at $q_0\in S\setminus C$ and,
for any $p\in\pi^{-1}(q_0)$, let
$$H(p):=\big\{g\in\PU V\mid gp=ph\text{ for some }h\in
H(q_0)\big\}\le\PU V$$
stand for the corresponding subgroup. (By convention, $H(q_0)$ acts on
$\pi^{-1}(q_0)$ from the right and this action commutes with that of
$\PU V$.)

Consider a piecewise vertical-horizontal path $c:[a,b]\to S$ and its
lift $\widetilde c:[a,b]\to E$, $\pi\circ\widetilde c=c$, that begins
at an arbitrary $p\in\pi^{-1}\big(c(a)\big)$, $\widetilde c(a)=p$. By
construction, the path $\widetilde c$ lies entirely in some $E_F$.
Denote by $\varphi_c:[a,b]\to\PGL V$ the path given by
$p\varphi_c(s)=\widetilde c(s)$ in terms of representatives
$p_1,p_2,p_3\in V$ with standard Gram matrix such that
$p=(p_1,p_2,p_3)$, where
$(p_1,p_2,p_3)\varphi_c(s):=\big(\varphi_c(s)^{-1}p_1,
\varphi_c(s)^{-1}p_2,\varphi_c(s)^{-1}p_3\big)$.
As it is easy to see, the lift of $c$ that begins at another point
$gp\in\pi^{-1}\big(c(a)\big)$, $g\in\PU V$, provides the same
$\varphi_c$.

Let $c':[b,d]\to S$ be one more piecewise vertical-horizontal path.
Then, of course,
$\varphi_{c\cup c'}(d)=\varphi_c(b)\circ\varphi_{c'}(d)$. In this way,
we can explicitly find the paths $\varphi_c$'s for piecewise
vertical-horizontal $c$'s if we know $\varphi_c$ in the simple cases of
a vertical/horizontal $c$.

Pick an arbitrary {\it basic\/} point $q_0\in S$ and consider closed
piecewise vertical-horizontal paths $c$ that begin at $q_0$. The
corresponding $\varphi_c$'s form a subgroup $H_0(q_0)\le\PGL V$ called
{\it rectangle holonomy group\/} at $q_0$. For any $p\in\pi^{-1}(q_0)$,
denote
$H_0(p):=\big\{g\in\PU V\mid gp=ph\text{ for some }h\in
H_0(q_0)\big\}\le\PU V$.
Since the lift of any vertical-horizontal path lies entirely in some
$E_F$, we obtain $H_0(p)\le C(F_p)$.

The rectangle holonomy group $H_0(q'_0)$ at another basic point
$q'_0\in S$ is isomorphic to $H_0(q_0)$ by means of the conjugation by
$\varphi_c$, where $c$ is a piecewise vertical-horizontal path from
$q_0$ to $q'_0$, $H_0(q_0)=\varphi_cH_0(q'_0)\varphi_c^{-1}$. As
$H_0(q_0)\simeq H_0(p)\le C(F_p)$ is commutative by Remark 3.9, this
isomorphism is independent of the choice of $c$.

We can contract any closed piecewise vertical-horizontal path based at
$q_0$, keeping it piecewise vertical-horizontal and based at $q_0$
during the contraction. Therefore, $H_0(q_0)$ is path-connected. Any
piecewise smooth closed path based at $q_0\in S\setminus C$ can be
approximated by a piecewise vertical-horizontal one. So,~$H_0(q_0)$ is
dense in $H(q_0)$. Since $H(q_0)$ is known to be a connected Lie group,
we conclude that $H(p)$ is a closed connected subgroup in $C(F_p)$ for
any $p\in\pi^{-1}(q_0)$. Moreover, $H(q_0)=H_0(q_0)$. Indeed,
as~$H(q_0)$ is abelian and $H_0(q_0)$ is path-connected, the subgroup
$H_0(q_0)$ is the image of a linear subspace under the exponential map.
It remains to observe that a small piecewise smooth closed path based
at $q_0\in S\setminus C$ can be approximated by a small piecewise
vertical-horizontal one.

Let $p$ be a real strongly regular triple and let $R$ denote the
corresponding real plane. Then $H_0(p)\le\Stab R$ because the bendings
involving $p_1,p_2$/$p_2,p_3$ keep the triple inside $R$. Since
$H_0(p)$ is path-connected, $H_0(p)\le\Stab^+R$, where
$\Stab^+R\le\PU V$ denotes the group of orientation-preserving
isometries of the hyperbolic plane $R$. So,
$H_0(p)\le\Stab^+R\cap C(F_p)$. It is immediate that
$\Stab^+R\cap C(F_p)$ is the $1$-parameter subgroup generated by $F_p$
in the group of isometries of the hyperbolic plane $R$.

Summarizing, we arrive at the following

\medskip

{\bf4.2.1.~Lemma.} {\sl The rectangle holonomy group coincides with the
holonomy group, i.e., $H_0(p)=H(p)$ for any strongly regular triple\/
$p\in E$. Moreover, $H_0(p)$ is a connected Lie subgroup in\/ $C(F_p)$,
$H_0(p)\le C(F_p)$, where\/ $F_p:=R(p_3)R(p_2)R(p_1)$.

If\/ $p$ is real, then $H_0(p)$ lies in the\/ $1$-parameter subgroup
generated by\/ $F_p$ in the group of isometries of the hyperbolic
plane\/ $R$, $H_0(p)\le\Stab^+R\cap C(F_p)$, where\/ $R$ stands for the
real plane spanned by\/ $p$}
$_\blacksquare$

\medskip

The main result of this paper is the following

\medskip

{\bf4.2.2.~Theorem.} {\sl Let\/ $p:=(p_1,p_2,p_3)$ be a strongly
regular triple and let\/ $F_p:=R(p_3)R(p_2)R(p_1)$. Then\/
$H_0(p)=C(F_p)$ unless the triple is real. If the triple is real,
then\/ $H_0(p)=\Stab^+R\cap C(F_p)$, where\/ $R$ stands for the real
plane of the triple.}

\medskip

{\bf4.2.3.~Tangent space to $E_F$.} In order to prove Theorem 4.2.2, we
are going to explicitly find the curvature tensor of the connection on
$E_F\to S$ and, using the Ambrose-Singer theorem, to calculate the
holonomy group $H(p)$. So, we fix $F$.

A tangent vector $t\in\T_pE_F$ can be described as a triple
$t:=(t_1,t_2,t_3)$ such that
$$\textstyle\frac
d{d\varepsilon}\big|_{\varepsilon=0}R\big(p_3+\varepsilon
t_3(p_3)\big)R\big(p_2+\varepsilon t_2(p_2)\big)R\big(p_1+\varepsilon
t_1(p_1)\big)=0,\leqno{\bold{(4.2.4)}}$$
where $p=(p_1,p_2,p_3)$ and
$t_i\in\T_{p_i}\Bbb P_\Bbb CV\subset\Lin_\Bbb C(V,V)$. One can easily
verify that
$$\textstyle\frac
d{d\varepsilon}\big|_{\varepsilon=0}R\big(p_i+\varepsilon
t_i(p_i)\big)=2(t_i+t_i^*)=2\widehat t_iR(p_i)\in{\T}_{R(p_i)}\SU V=\su
V\cdot R(p_i)$$
and that $\widehat t_jR(p_j)+R(p_j)\widehat t_j=0$ (see Definition 2.3
for $\widehat t$). Hence, (4.2.4) is equivalent to
$$\widehat t_3R(p_3)R(p_2)R(p_1)+R(p_3)\widehat
t_2R(p_2)R(p_1)+R(p_3)R(p_2)\widehat t_1R(p_1)=0,$$
i.e., to
$$-\widehat t_3+\widehat t_2+R(p_2)\widehat
t_1R(p_2)=0.\leqno{\bold{(4.2.5)}}$$
We have reduced the task of describing the tangent space $\T_pE_F$ to
the following one. Given a regular triple $p'=(p'_2,p'_1,p'_3)$ (where
$p'_1:=R(p_2)p_1$, $p'_2=p_2$, and $p'_3=p_3$; note that
$F=R(p'_3)R(p'_1)R(p'_2)$), we look for $t'_j\in\T_{p'_j}\Bbb P\Bbb CV$
(where $t'_1:=R(p_2)t_1$, $t'_2:=t_2$, and $t'_3:=-t_3$) such that
$\widehat{t'}_3+\widehat{t'}_2+\widehat{t'}_1=0$.

Let $[g'_{jk}]$ denote the Gram matrix of the $p'_j$'s and let
$t'_j:=\sum_k\langle-,p'_j\rangle a_{jk}p'_k$. The equation (4.2.5) is
equivalent to the equalities $a_{kj}=\overline a_{jk}$ and the
conditions $\sum_ka_{jk}p'_k\in{p'}_j^\perp$ are equivalent to the
equalities $\sum_ka_{jk}g'_{kj}=0$. They can be rewritten as
$a_{jj}:=\frac1{g'_{jj}}\Re\big(\sum_{k\ne j}a_{jk}\big)$ and
$\Im\big(\sum_{k\ne j}a_{jk}g'_{kj}\big)=0$ for all $j$. The last three
equations in the $a_{jk}$'s, $j\ne k$, are linearly dependent (their
sum equals zero) and define a $4$-dimensional space of solutions
because $g'_{jk}\ne0$.

\medskip

{\bf4.2.6.~Lemma.} {\sl The vector fields
$$b_1(p):=\big(\langle-,p_1\rangle(\textstyle
\frac{p_1}{g_{11}}-\frac{p_2}{g_{21}}),\langle-,p_2\rangle
(\frac{p_1}{g_{12}}-\frac{p_2}{g_{22}}),0\big),\qquad
b_2(p):=\big(0,\langle-,p_2\rangle(\frac{p_2}{g_{22}}-
\frac{p_3}{g_{32}}),\langle-,p_3\rangle(\frac{p_2}{g_{23}}-
\frac{p_3}{g_{33}})\big)$$
defined for all\/ $p:=(p_1,p_2,p_3)\in E$ never vanish and span the
horizontal\/ $\Bbb R$-linear space of the connection on\/ $E\to S$
over\/ $S\setminus C$, where\/ $[g_{jk}]$ stands for the Gram matrix of
the\/ $p_j$'s. Denote\/ $\tau:=\frac{g_{13}g_{22}}{g_{12}g_{23}}$. Then
$$[b_1,b_2](p)=\Big(\langle-,p_1\rangle
\overline\tau(\textstyle\frac{p_2}{g_{21}}-\frac{p_3}{g_{31}}),
\langle-,p_2\rangle\big((2-\tau)\frac{p_1}{g_{12}}+
2i\Im\tau\frac{p_2}{g_{22}}+(\overline\tau-2)\frac{p_3}{g_{32}}\big),
\langle-,p_3\rangle\tau(\frac{p_1}{g_{13}}-\frac{p_2}{g_{23}})\Big).$$}

{\bf Proof.} It is immediate that the vectors $b_1$ and $b_2$ never
vanish and are tangent to the lifts of the vertical and horizontal
lines, respectively.

Let $f_j$ be an analytic function on $E$ depending only on
$p_j\in\Bbb P_\Bbb CV$, $j=1,2,3$. Then
$$b_1f_1=\textstyle\frac d{d\varepsilon}\big|_{\varepsilon=0}
f_1(p_1-\varepsilon\frac{g_{11}}{g_{21}}p_2),\qquad b_1f_2=\frac
d{d\varepsilon}\big|_{\varepsilon=0}f_2(p_2+
\varepsilon\frac{g_{22}}{g_{12}}p_1),\qquad b_1f_3=0,$$
$$b_2f_1=0,\qquad b_2f_2=\textstyle\frac
d{d\varepsilon}\big|_{\varepsilon=0}f_2(p_2-
\varepsilon\frac{g_{22}}{g_{32}}p_3),\qquad b_2f_3=\frac
d{d\varepsilon}\big|_{\varepsilon=0}
f_3(p_3+\varepsilon\frac{g_{33}}{g_{23}}p_2),$$
$$[b_1,b_2]f_1=-b_2b_1f_1=-\textstyle\frac
d{d\delta}\big|_{\delta=0}\frac
d{d\varepsilon}\big|_{\varepsilon=0}f_1\Big(p_1-\varepsilon
g_{11}\displaystyle\frac{\frac{p_2}{g_{21}}+\delta\frac{g_{22}}{g_{21}}
(\frac{p_2}{g_{22}}-\frac{p_3}{g_{32}})}{1+\delta\frac{g_{22}}{g_{21}}
(\frac{g_{21}}{g_{22}}-\frac{g_{31}}{g_{32}})}\Big)=
\langle-,p_1\rangle\textstyle\frac{g_{22}g_{31}}{g_{32}g_{21}}
(\frac{p_2}{g_{21}}-\frac{p_3}{g_{31}})f_1.$$
Similarly,
$[b_1,b_2]f_3=
\langle-,p_3\rangle\tau(\frac{p_1}{g_{13}}-\frac{p_2}{g_{23}})f_3$.

We have calculated the first and the last components of the vector
$[b_1,b_2](p)$. Since the vectors $b_1(p'),b_2(p')$ are tangent to
$E_{F_p}$ for all $p'\in E_{F_p}$, the vector $[b_1,b_2](p)$ has to be
tangent to~$E_{F_p}$. It~is straightforward that the right-hand side of
the expression given in Lemma 4.2.6 satisfies the equation~(4.2.5),
which uniquely determines the second component
$_\blacksquare$

\medskip

{\bf4.2.7.~Curvature tensor of $E_F\to S$.} Let $l\in\su V$. By [AGG,
Lemma 4.1.4], the associated to $l$ tangent vector to the fibre of
$\pi:E\to S$ at $p:=(p_1,p_2,p_3)\in E$ has the form
$$\big(\langle-,p_1\rangle\textstyle\frac{\pi[p_1]l(p_1)}{g_{11}},
\langle-,p_2\rangle\frac{\pi[p_2]l(p_2)}{g_{22}},
\langle-,p_3\rangle\frac{\pi[p_3]l(p_3)}{g_{33}}\big),$$
where $\pi[p_j]:V\to p_j^\perp$ denotes the orthogonal projection and
$[g_{jk}]$ is the Gram matrix of the $p_j$'s. Conversely, let
$t:=\big(\langle-,p_1\rangle q_1,\langle-,p_2\rangle
q_2,\langle-,p_3\rangle q_3\big)$
be a tangent vector to $E$ at $p:=(p_1,p_2,p_3)\in E$ with
$q_j\in p_j^\perp$. Denote by $T$ the $\Bbb C$-linear map $T:V\to V$
given by $T(p_j):=g_{jj}q_j$ and consider a $\Bbb C$-linear map
$L:V\to V$ such that $L(p_j):=T(p_j)-id_jp_j$ with $d_j\in\Bbb R$ for
all $j=1,2,3$. Then $t$ is tangent to the fibre iff $L\in\su V$ for
suitable $d_j$'s. In terms of the basis $p_1,p_2,p_3$, this means that
$\tr T=i\tr D$ and
$T^t[g_{jk}]+[g_{jk}]\overline T=i\big(D[g_{jk}]-[g_{jk}]D\big)$, where
$D:=\left[\smallmatrix
d_1&0&0\\0&d_2&0\\0&0&d_3\endsmallmatrix\right]$.
In other words, $t$ is tangent to the fibre of $\pi:E\to S$ iff
$\tr T\in i\Bbb R$ and
$$T^t[g_{jk}]+[g_{jk}]\overline
T=i\left[\smallmatrix0&a_{12}g_{12}&a_{13}g_{13}\\a_{21}g_{21}&0&
a_{23}g_{23}\\a_{31}g_{31}&a_{32}g_{32}&0\endsmallmatrix\right]
\leqno{\bold{(4.2.8)}}$$
with $a_{kj}=-a_{jk}\in\Bbb R$ and $a_{12}+a_{23}+a_{31}=0$. Moreover,
in this case, $d_1=\frac{d+a_{12}-a_{31}}3$,
$d_2=\frac{d+a_{23}-a_{12}}3$, $d_3=\frac{d+a_{31}-a_{23}}3$, where
$id:=\tr T$. In particular,
$$d_1=\textstyle\frac{d+2a_{12}+a_{23}}3.\leqno{\bold{(4.2.9)}}$$

Denote by $\omega$ and $\Omega$ the $1$-form of the connection on
$E_F\to S$ and its curvature tensor, both defined over $S\setminus C$.
They take values in the Lie algebra $c(F)\le\su V$ of the group $C(F)$.
Since $\Omega(b_2\wedge b_1)=\omega[b_1,b_2]$ by Lemma 4.2.6, we need
the following lemma.

\medskip

{\bf4.2.10.~Lemma.} {\sl Let\/
$p:=(p_1,p_2,p_3)\in E_F\setminus\pi^{-1}(C)$. Then
$$\omega[b_1,b_2](p_1)=g_{11}\textstyle\frac{(1-\beta-t_2)t_1t_2-
2\alpha^2}{t_1^2t_2^2(t-1)}(\frac{p_1}{g_{11}}-\frac{p_2}{g_{21}})+
g_{11}\overline\tau(\frac{p_2}{g_{21}}-\frac{p_3}{g_{31}})
+\frac{i\alpha(1-\beta-3t_2+2t_1t_2)}{3t_1^2t_2^2(t-1)}p_1,$$
where\/ $\tau:=\frac{g_{13}g_{22}}{g_{12}g_{23}}$ and\/ $[g_{jk}]$
stands for the Gram matrix of the\/ $p_j$'s.}

\medskip

{\bf Proof.} By Lemma 4.2.6, there exist unique $c_1,c_2\in\Bbb R$ such
that the vector $c_1b_1(p)+c_2b_2(p)+[b_1,b_2](p)$ is tangent to the
fibre $E_F\cap\pi^{-1}(q)$ of $E_F\to S$. By Lemma 4.2.6 and the
considerations in 4.2.7,
$\omega[b_1,b_2](p_1)=g_{11}c_1(\frac{p_1}{g_{11}}-\frac{p_2}{g_{21}})+
g_{11}\overline\tau(\frac{p_2}{g_{21}}-\frac{p_3}{g_{31}})-id_1p_1$.

The matrices of $b_1(p),b_2(p),[b_1,b_2](p)$ are
$$B_1:=\left[\smallmatrix1&\frac{g_{22}}{g_{12}}&0\\
-\frac{g_{11}}{g_{21}}&-1&0\\0&0&0\endsmallmatrix\right],\qquad
B_2:=\left[\smallmatrix0&0&0\\0&1&\frac{g_{33}}{g_{23}}\\
0&-\frac{g_{22}}{g_{32}}&-1\endsmallmatrix\right],\qquad
B:=\left[\smallmatrix0&\frac{g_{22}}{g_{12}}(2-\tau)&
\frac{g_{33}}{g_{13}}\tau\\\frac{g_{11}}{g_{21}}\overline\tau&
2i\Im\tau&-\frac{g_{33}}{g_{23}}\tau\\-\frac{g_{11}}
{g_{31}}\overline\tau&\frac{g_{22}}{g_{32}}(\overline\tau-2)
&0\endsmallmatrix\right].$$
We calculate the hermitian matrices
$$B_1^t[g_{jk}]+[g_{jk}]\overline B_1=\left[\smallmatrix0&0&*\\
0&0&g_{23}(\tau-1)\\g_{31}(1-\frac{g_{11}g_{32}}{g_{31}g_{12}})&*&0
\endsmallmatrix\right],\qquad B_2^t[g_{jk}]+[g_{jk}]\overline
B_2=\left[\smallmatrix0&g_{12}(1-\tau)&*\\*&0&0\\g_{31}(\frac{g_{21}
g_{33}}{g_{23}g_{31}}-1)&0&0\endsmallmatrix\right],$$
$$B^t[g_{jk}]+[g_{jk}]\overline B=\left[\smallmatrix0&g_{12}
\big(\frac{g_{11}g_{22}}{g_{12}g_{21}}+\tau(\tau-3)+\overline\tau\big)&
*\\*&0&-g_{23}\big(\frac{g_{22}g_{33}}{g_{23}g_{32}}+\tau(\tau-3)+
\overline\tau\big)\\g_{31}\frac{g_{13}g_{22}}{g_{12}g_{23}}
(\frac{g_{11}g_{32}}{g_{31}g_{12}}-\frac{g_{21}g_{33}}{g_{23}g_{31}})&
*&0\endsmallmatrix\right].$$

Let $T:=c_1B_1+c_2B_2+B$. Then $\tr T=id$, where $d:=2t'$ and
$t+it':=\tau$.

The $12$- and $23$-components of the matrix
$T^t[g_{jk}]+[g_{jk}]\overline T$ belong to $ig_{12}\Bbb R$ and
$ig_{12}\Bbb R$ iff
$$c_2(t-1)=\Re(\textstyle\frac{g_{11}g_{22}}{g_{12}g_{21}}+\tau(\tau-3)
+\overline\tau),\qquad c_1(t-1)=\Re(\frac{g_{22}g_{33}}{g_{23}g_{32}}+
\tau(\tau-3)+\overline\tau),$$
respectively. Moreover, we can also see that
$$a_{12}=\Im\big(-c_2\tau+\tau(\tau-3)+\overline\tau\big),\qquad
a_{23}=\Im\big(c_1\tau-\tau(\tau-3)-\overline\tau\big).$$
Hence,
$$c_1=\textstyle\frac{\frac1{t_2}+t^2-{t'}^2-2t}{t-1},\qquad
c_2=\frac{\frac1{t_1}+t^2-{t'}^2-2t}{t-1},\qquad
a_{12}=-c_2t'+2tt'-4t',\qquad a_{23}=c_1t'-2tt'+4t'.$$
By (4.2.9),
$d_1=\frac{(c_1-2c_2+2t-2)t'}3=\frac{(\frac1{t_2}-\frac2{t_1}+t^2+
{t'}^2-2t+2)t'}{3(t-1)}$.

From
$t_1t_2\overline\tau=\frac{g_{12}g_{23}g_{31}}{g_{11}g_{22}g_{33}}$, we
conclude that $t'=-\frac\alpha{t_1t_2}$. Using equation (4.1.2), we can
see that $c_1=\frac{(1-\beta-t_2)t_1t_2-2\alpha^2}{t_1^2t_2^2(t-1)}$
and $d_1=-\frac{\alpha(1-\beta-3t_2+2t_1t_2)}{3t_1^2t_2^2(t-1)}$
$_\blacksquare$

\medskip

{\bf Proof of Theorem 4.2.2.} By Lemma 4.2.10, the holonomy group is
not trivial. So, we assume that $\alpha\ne0$.

We fix the point $p_1$ and the geodesic $\G{\wr}p_2,p_3{\wr}$. It
suffices to show that, for a generic $(p_1,p_2,p_3)\in E_F$, the
element $\omega[b_1,b_2](p_1)\in V$ does not remain
$\Bbb R$-proportional when we bend $p_2,p_3$.

Let $\{u_3\}:=\Bbb P_\Bbb Cp_1^\perp\cap\L(p_2,p_3)$ and let $u_2$ be
the polar point to $\L(p_1,u_3)$. Then
$\frac{\langle p_2,u_2\rangle}{g_{21}}=\frac{\langle
p_3,u_2\rangle}{g_{31}}$
because
$\frac{p_2}{g_{21}}-\frac{p_3}{g_{31}}\in\Bbb P_\Bbb
Cp_1^\perp\cap\L(p_2,p_3)$.
Note that $b:=\frac{\langle g,u_2\rangle}{\langle g,p_1\rangle}$ is
independent of the choice of
$g\in\G{\wr}p_2,p_3{\wr}\setminus\Bbb P_\Bbb Cp_1^\perp$. Indeed,
$\frac{\langle g,u_2\rangle}{\langle
g,p_1\rangle}=\frac{(1+\frac{rg_{31}}{g_{32}g_{21}})\langle
p_2,u_2\rangle}{g_{21}+\frac{rg_{31}}{g_{32}}}=\frac{\langle
p_2,u_2\rangle}{g_{21}}$
because
$\langle p_3,u_2\rangle=\frac{g_{31}\langle p_2,u_2\rangle}{g_{21}}$,
where $g:=p_2+\frac{rp_3}{g_{32}}\in\G{\wr}p_2,p_3{\wr}$ for any
$r\in\Bbb R$. The point $u_2$ cannot be polar to $\L(p_2,p_3)$.
Therefore, $b\ne0$.

It remains to observe that
$$\big\langle\omega[b_1,b_2](p_1),u_2\big\rangle=\textstyle
\frac{2\alpha^2-(1-\beta-t_2)t_1t_2}{t_1^2t_2^2(t-1)}g_{11}b,\qquad
\big\langle\omega[b_1,b_2](p_1),p_1\big\rangle=
\frac{i\alpha(1-\beta-3t_2+2t_1t_2)}{3t_1^2t_2^2(t-1)}g_{11}$$
by Lemma 4.2.10
$_\blacksquare$

\medskip

{\bf4.3.~Regular triples.} In this subsection, we discuss regular
triples that are not strongly regular~and, a little bit, spherical
configurations.

Consider a regular triple $p_1,p_2,p_3\in\L$ that lies in a projective
line $\L$. The case when all points are negative is in fact covered by
Theorem 4.2.2 : it does not differ from the case of a real strongly
regular triple because the geometries of the Poincar\'e disc and of the
Beltrami-Klein disc are essentially the same. When a point $p_j$ is
positive,  we can substitute  $p_j$ by the negative point $p'_j\in\L$
orthogonal to $p_j$. Indeed, let $p$ denote the point polar to $\L$.
Then $R(p)$ commutes with every $R(p_k)$ and $R(p)R(p_j)=R(p'_j)$.
Every geodesic that passes through $p_j$ passes necessarily through
$p'_j$ [AGr, Section 3]. Now, it is easy to see that a bending
involving $p_j,p_k$ corresponds to a bending involving $p'_j,p_k$.

Let $p_1,p_2,p_3\in R$ be a real regular triple with a positive $p_j$,
where $R$ stands for the real plane spanned by the $p_k$'s. In terms of
the Beltrami-Klein disc $D$ of all negative points in $R$, the set
$\G:=D\cap\Bbb P_\Bbb Cp_j^\perp$ is a geodesic. Denote by $\Gamma$ the
geodesic in $D$ that passes through $p_k$ and is perpendicular to $\G$.
When we bend $p_j,p_k$, we in fact move $\G$ and $p_k$ keeping the
distance between $\G$ and $p_k$, keeping $\G$ orthogonal to~$\Gamma$,
and keeping $p_k\in\Gamma$. At some moment, the third point $p_l$ lies
in $\G$, $p_l\in\G$. If necessary, we can replace $R(p_j)$ and $R(p_k)$
by $R\big(R(p_j)p_k\big)$ and $R(p_j)$. Now, two subsequent points
$p_j,p_l$ are orthogonal and $R(p_j)R(p_l)$ equals some $R(p)$ due to
orthogonal relations.

If one discards the requirement that at most one point is positive,
there is an example of a triple $p_1,p_2,p_3$ such that the projective
lines $\L(p_1,p_2),\L(p_2,p_3)$ are hyperbolic and $\L(p_2,p_3)$
becomes spherical after a bending involving $p_1,p_2$. Indeed, there is
a number $z\in\Bbb C$ such that $|\frac14+z|<\frac12$ and $1<|1+z|$
(for instance, $z=\frac18$). By Sylvester's Criterion, there are points
$v_1,v_2,p_3\in V$ with the Gram matrix
$\left[\smallmatrix0&\frac12&1\\\frac12&0&\overline z\\1&z&1
\endsmallmatrix\right]$.
We put
$$p_1:=2v_1-\textstyle\frac12v_2,\qquad p_2:=v_1+v_2,\qquad
p'_1:=v_1-v_2,\qquad p'_2:=\frac12v_1+2v_2.$$
Then the projective lines $\L(p_1,p_2),\L(p_2,p_3)$ are hyperbolic
because $1<|1+z|$, whereas the projective line $\L(p'_2,p_3)$ is
spherical because $|\frac12+2z|<1$.

Nevertheless, in principle, there can exist a {\it spherical\/}
configuration of points $p_1,p_2,\dots,p_n$, i.e., a~nontrivial
configuration subject to the relation $R(p_n)\dots R(p_2)R(p_1)=1$ with
spherical projective lines $\L(p_{j-1},p_j)$ (the indices are modulo
$n$) that remain spherical after a finite number of bendings. In this
case, one should suitably modify Conjecture 1.1.

While looking for basic relations, one should not go too far. Most
likely, every isometry in $\SU V$ is a product of $5$ reflections.
Hence, it must be sufficient to deal with $n\le11$. However, it is
quite possible that actually we only need to study the relations with
$n\le6$.

\bigskip

\centerline{\bf5.~Pentagons}

\medskip

In this section, we begin to study $3$ basic relations (in the sense of
Conjecture 1.1), the pentagons. (Yet, we do not know if there exists a
spherical configuration with $n=5$.) A {\it pentagon\/} is a
configuration $p_1,p_2,p_3,p_4,p_5\in\Bbb P_\Bbb CV\setminus\S V$ of
$5$ nonisotropic points such that at most one of the $p_j$'s is
positive, $p_{j-1}$~is not orthogonal nor equal to $p_j$ for all $j$
(the indices are modulo $5$), and
$$R(p_5)R(p_4)R(p_3)R(p_2)R(p_1)=\delta$$
in $\SU V$, where $\delta^3=1$. We will show that pentagons considered
modulo congruence form exactly $3$ components. Moreover, two pentagons
with the same $\delta$ are connected by means of a finite number of
bendings. Each component is a smooth $4$-manifold.

If $\delta=1$, the pentagons are known to provide faithful and discrete
representations $\varrho:H_5\to\PU V$ [ABG], where $H_n$ denotes the
{\it hyperelliptic\/} group, i.e., the group generated by
$r_1,r_2,\dots,r_n$ with the defining relations $r_j^2=1$,
$j=1,2,\dots,n$, and $r_n\dots r_2r_1=1$. Moreover, in this case the
$p_j$'s are all negative and span a real plane $R$. In other words, the
representation $\varrho$ is $\Bbb R$-fuchsian as $R$ is
$\varrho H_5$-stable. So, we arrive at a sort of Toledo rigidity [Tol]
: while deforming a pentagon with $\delta=1$, the~representation
$\varrho$ remains faithful, discrete, and $\Bbb R$-fuchsian.

The two components with $\delta\ne1$ are congruent: one can be obtained
from the other by reflecting the pentagons in a real plane (they are
`complex conjugated'). One of the $p_j$'s is always positive. In [Ana],
we construct a pentagon with $\delta\ne1$ that provides a discrete and
faithful representation $\varrho:H_5\to\PU V$. We believe that every
pentagon provides a faithful and discrete representation (Conjecture
1.2).

\medskip

{\bf5.1.~Signs of points.} It follows from
$R(p_3)R(p_2)R(p_1)=\delta R(p_4)R(p_5)$ and Lemma 3.1 that
$$8i\alpha(p_1,p_2,p_3)+4\beta(p_1,p_2,p_3)-1=
\delta\big(4\ta(p_4,p_5)-1\big).\leqno{\bold{(5.2)}}$$

Suppose that $\delta=1$ and $p_1,p_2,p_3\in\B V$. Then
$\alpha(p_1,p_2,p_3)=0$, $\beta(p_1,p_2,p_3)\ge0$, and
$\beta(p_1,p_2,p_3)=\ta(p_4,p_5)$, implying that $p_4,p_5$ have the
same sign, i.e., are negative.

Suppose that $\delta\ne1$ and $p_4,p_5\in\B V$. Then
$\Re\delta=-\frac12$ and
$4\beta(p_1,p_2,p_3)-1=-\frac12\big(4\ta(p_4,p_5)-1\big)$. Hence,
$8\beta(p_1,p_2,p_3)=3-4\ta(p_4,p_5)$. Since $\ta(p_4,p_5)>1$, we
obtain $\beta(p_1,p_2,p_3)<0$. So, one of the points $p_1,p_2,p_3$ is
positive.

\medskip

{\bf5.3.~Existence.} For distinct $p_4,p_5\in\B V$, the eigenvectors of
$R(p_4)R(p_5)$ are the vertices of the geodesic $\G{\wr}p_4,p_5{\wr}$
and the point polar to $\L(p_4,p_5)$. The corresponding eigenvalues are
$s^{-1},s,1$, where $s>\nomathbreak1$. By Proposition 3.8, we can find
a regular triple $p_1,p_2,p_3$ such that
$R(p_3)R(p_2)R(p_1)=\delta R(p_4)R(p_5)$. If~$p_3=p_4$ or $p_1=p_5$,
then $p_1,p_2$ or $p_2,p_3$ are orthogonal by Remark 3.2, which
contradicts the fact that the triple $p_1,p_2,p_3$ is regular. If
$p_1,p_5$ are orthogonal, then, after applying the orthogonal relation,
we obtain $R(p_4)R(p_3)R(p_2)R(p)=\delta$ with $p,p_1,p_5$ pairwise
orthogonal. By Corollary 3.3, $p,p_2,p_3,p_4$ lie on the geodesic
$\G{\wr}p_2,p_3{\wr}$. As $p,p_1$ are orthogonal,
$p_1\in\G{\wr}p_2,p_3{\wr}$, which contradicts the fact that the triple
$p_1,p_2,p_3$ is regular. Symmetric arguments work if $p_3,p_4$ are
orthogonal.

Thus, we have shown that, for any $\delta$ with $\delta^3=1$, there
exists a pentagon.

\medskip

{\bf5.4.~Connectedness by means of bendings.} Given a pentagon, the
triple $p_1,p_2,p_3$ is strongly regular. Indeed, if $p_1,p_2,p_3$ are
on a same geodesic $\G$, then $R(p_3)R(p_2)R(p_1)=R(p)$ with $p\in\G$
by Corollary 3.3. Now, $R(p_5)R(p_4)R(p)=\delta$ implies that $p_4,p_5$
are orthogonal by Remark 3.2. A~contradiction. If $p_1,p_2,p_3$ are on
a same projective line, then $\beta(p_1,p_2,p_3)=0$ and (5.2) implies
that $\Re\delta=-\frac12$ and, hence, $4\ta(p_4,p_5)=3$. So, the
projective line $\L(p_4,p_5)$ is spherical; a contradiction. If
$p_1,p_2,p_3$ are on a same real plane, then $\alpha(p_1,p_2,p_3)=0$.
Therefore, $\delta=1$ by (5.2) and the $p_j$'s are all negative by 5.1.

In particular, after any bending, a pentagon remains a pentagon.
Indeed, after a bending involving, say, $p_2,p_3$, the triples
$p_1,p_2,p_3$ and $p_2,p_3,p_4$ remain strongly regular by Corollary
4.1.5.

Given two pentagons with the same $\delta$ and negative
$p_2,p_3,p_4,p_5$, by means of a bending involving $p_3,p_4$, one can
make $\ta(p_4,p_5)$ arbitrarily big. So, we can assume that
$\ta(p_4,p_5)$ is the same for both pentagons. The connectedness
follows now from Corollary 4.1.5.

\medskip

{\bf5.5.~Real pentagons.} Consider the case of $\delta=1$. By 5.1 and
5.4, the points $p_1,p_2,p_3$ form a strongly regular real triple and,
therefore, span a real plane $R$, which is stable under the isometry
$I:=R(p_3)R(p_2)R(p_1)=R(p_4)R(p_5)$. As the fixed points of the
isometry $R(p_4)R(p_5)$ are the (isotropic) vertices $v,v'$ of the
geodesic $\G{\wr}p_4,p_5{\wr}$ and the (positive) point polar to the
projective line $\L(p_4,p_5)$, the~isometry $I$ of the Beltrami-Klein
disc $D\subset R$ has to be hyperbolic. So, its isotropic fixed points
should coincide with $v,v'$, implying $v,v'\in R$ and, hence,
$p_4,p_5\in R$. As there is no essential difference between the
Beltrami-Klein disc and the Poincar\'e one, we can apply [ABG,
Corollary 3.16] and conclude that the pentagon provides a faithful and
discrete representation $\varrho:H_5\to\PU V$.

\medskip

{\bf5.6.~Smoothness.} The space of real pentagons is smooth and
$4$-dimensional by [ABG, Corollary~3.17]. So, we assume
$\delta:=-\frac12+i\frac{\sqrt3}2$ and $p_1$ positive.

Let $t_j:=\ta(p_j,p_{j+1})$ (the indices are modulo $5$). It follows
from (5.2) that $\alpha(p_1,p_2,p_3)=\frac{\sqrt3(4t_4-1)}{16}$ and
$\beta(p_1,p_2,p_3)=\frac{3-4t_4}8$. The inequalities (4.1.3) for
$t_1,t_2,t_4$ are
$$t_1<0,\qquad1<t_2,\qquad1<t_4.\leqno{\bold{(5.7)}}$$
Therefore, $\beta<0$, as desired. The equation (4.1.2) takes the form
$$(t_1-1)(t_2-1)=t_1t_2(t-1)^2+\textstyle\frac{3(4t_4-1)^2}
{256t_1t_2}+\frac{3-4t_4}8.\leqno{\bold{(5.8)}}$$
The equation (5.8) in $t_1,t_2,t_4,t$ together with the inequalities
(5.7) defines a smooth $3$-manifold $T\subset\Bbb R^4(t_1,t_2,t_4,t)$
fibred by Lemma 4.1.1 into planes over the axis $A:=\{t_4\mid t_4>1\}$.
Hence, $T\simeq\Bbb R^3$. Geometrically, a pentagon is determined by a
point in $T$ uniquely up to a bending involving $p_4,p_5$. In~other
words, the space of pentagons with a given $\delta$ is diffeomorphic to
$\Bbb R^4$.

\newpage

\centerline{\bf6.~References}

\medskip

[ABG] S.~Anan$'$in, E.~C.~B.~Gon\c calves, {\it A hyperelliptic view on
Teichm\"uller space.~{\rm I},} preprint\newline
http://arxiv.org/abs/0709.1711

[AGG] S.~Anan$'$in, C.~H.~Grossi, N.~Gusevskii, {\it Complex hyperbolic
structures on disc bundles over surfaces,}
Int.~Math.~Res.~Not.~{\bf2011} (2011), no.~19, 4295--4375, see also
http://arxiv.org/abs/math/0511741

[AGr] S.~Anan$'$in, C.~H.~Grossi, {\it Coordinate-free classic
geometries,} Mosc.~Math.~J.~{\bf11} (2011), 633--655, see also
http://arxiv.org/abs/math/0702714

[Ana] S.~Anan$'$in, {\it A discrete pentagon,} in preparation

[Gol1] W.~M.~Goldman, {\it Topological components of spaces of
representations,} Invent.~Math.~{\bf93} (1988), no.~3, 557--607

[Gol2] W.~M.~Goldman, {\it Complex hyperbolic geometry,} Oxford
Mathematical Monographs. Oxford Science Publications. The Clarendon
Press, Oxford University Press, New York, 1999, xx+316 pp.

[Pra] A.~Pratoussevitch, {\it Traces in complex hyperbolic triangle
groups,} Geometriae Dedicata {\bf111} (2005), 159--185

[Tol] D.~Toledo, {\it Representations of surface groups in complex
hyperbolic space,} J.~Differential Geom. {\bf29} (1989), no.~1,
125--133

\enddocument